\documentclass[final]{siamart171218}
\usepackage{amsmath,amssymb,stmaryrd, mathtools}
\usepackage{footmisc}
\usepackage{hyperref}
\usepackage{graphicx}
\usepackage{caption}
\usepackage{subcaption}
\usepackage{wrapfig}
\usepackage{xcolor}
\usepackage{cleveref}
\usepackage{accents}
\usepackage{cases}
\usepackage{enumerate}
\usepackage{comment}
\usepackage{tikz-cd}
\usepackage{array}

\usepackage[normalem]{ulem}

\newcommand{\ccrf}[2]{#1}
\newcommand{\ccbf}[2]{#1}

\newcommand{\bs}[1]{\boldsymbol{#1}}

\newcommand{\pp}[2]{\frac{\partial #1}{\partial #2}}

\newcommand{\vect}[1]{\mathbf{#1}}
\newtheorem{remark}{Remark}

\numberwithin{equation}{section}

\newcounter{savefootnote}
\newcounter{symfootnote}
\newcommand{\symfootnote}[1]{%
   \setcounter{savefootnote}{\value{footnote}}%
   \setcounter{footnote}{\value{symfootnote}}%
   \ifnum\value{footnote}>8\setcounter{footnote}{0}\fi%
   \let\oldthefootnote=\thefootnote%
   \renewcommand{\thefootnote}{\fnsymbol{footnote}}%
   \footnote{#1}%
   \let\thefootnote=\oldthefootnote%
   \setcounter{symfootnote}{\value{footnote}}%
   \setcounter{footnote}{\value{savefootnote}}%
}

\title{A Fast Algebraic Multigrid Solver and Accurate Discretization for Highly Anisotropic Heat Flux I: Open Field Lines \thanks{Corresponding author email address: gwimmer@lanl.gov
\funding{G. A. W., T. G., and X. T. have been supported by the U.S. Department of Energy Office of Fusion Energy Sciences and Office of Advanced Scientific Computing Research under the Tokamak Disruption Simulation (TDS) Scientific Discovery through Advanced Computing (SciDAC) project, as well as the Base Theory Program, both at Los Alamos National Laboratory (LANL) under contract No. 89233218CNA000001. B. S. S. was supported by the Laboratory Directed Research and Development program of Los Alamos National Laboratory as a Nicholas C. Metropolis Fellow and under project number 20220174ER. The computations have been performed using resources of the National Energy Research Scientific Computing Center (NERSC), a U.S. Department of Energy Office of Science User Facility operated under Contract No. DE-AC02-05CH11231. Los Alamos National Laboratory report number LA-UR-23-20885}}}

\author{Golo A. Wimmer\thanks{Theoretical Division, Los Alamos National Laboratory, USA.}
   \and Ben S. Southworth\footnotemark[2]
     \and Thomas J. Gregory\thanks{Imperial College London}
   \and Xian-Zhu Tang\footnotemark[2]
}

\begin{document}
\maketitle
\allowdisplaybreaks

\begin{abstract}
We present a novel solver technique for the anisotropic heat flux equation, aimed at the high level of anisotropy seen in magnetic confinement fusion plasmas. Such problems pose two major challenges: (i) discretization accuracy and (ii) efficient implicit linear solvers. We simultaneously address each of these challenges by constructing a new finite element discretization with excellent accuracy properties, tailored to a novel solver approach based on algebraic multigrid (AMG) methods designed for advective operators. \ccbf{We pose the problem in a mixed formulation, introducing the directional temperature gradient as an auxiliary variable. The temperature and auxiliary fields are discretized in a scalar discontinuous Galerkin space with upwinding principles used for discretizations of advection. We demonstrate the proposed discretization's superior accuracy over other discretizations of anisotropic heat flux, achieving error $1000\times$ smaller for anisotropy ratio of $10^9$, for \emph{closed field lines}. The block matrix system is reordered and solved in an approach where the two advection operators are inverted using AMG solvers based on approximate ideal restriction (AIR), which is particularly efficient for upwind discontinuous Galerkin discretizations of advection. To ensure that the advection operators are non-singular, in this paper we restrict ourselves to considering open (acyclic) magnetic field lines for the linear solvers. We demonstrate fast convergence of the proposed iterative solver in highly anisotropic regimes where other diffusion-based AMG methods fail.}{i}
\end{abstract}
\begin{keywords}
Anisotropic heat flux, anisotropic diffusion, auxiliary operator, algebraic multigrid, AIR
\end{keywords}



\section{Introduction}\label{sec:intro} The heat flux of plasmas in magnetic confinement fusion exhibits strong anisotropy \cite{jardin2010computational}. Heat can be conducted at a rate many orders of magnitude greater parallel to magnetic field lines than it is perpendicular to the magnetic field. \ccbf{For parallel and perpendicular conductivities $\kappa_\parallel$ and $\kappa_\perp$, respectively,}{1} this ratio $\frac{\kappa_{\parallel}}{\kappa_{\perp}}$ can exceed $10^{10}$ in realistic experiments \cite{gunter2007finite}. If a numerical simulation does not adequately resolve the anisotropy of the heat flux, the resulting simulated plasma confinement times may be spuriously short. Unfortunately, many numerical approaches suffer from the large parallel component of heat flux polluting the perpendicular direction, leading to non-physically large perpendicular heat flux. In addition, anisotropic heat flux in a magnetized plasma has a very restrictive CFL condition because of the large $\kappa_\parallel$, so for realistic applications an implicit solver is typically necessary. In the isotropic limit with $\kappa_\parallel = \kappa_\perp$, the problem reduces to the heat equation, and classical geometric and algebraic multigrid (AMG) solvers \cite{stuben1983algebraic} are known to be effective. Conversely, the highly anisotropic regime, $\kappa_\parallel \gg \kappa_\perp$, is notoriously difficult for multilevel and iterative solvers, and efficient parallel solvers for realistic applications remain a largely open question.

One common approach to address pollution of flux in the perpendicular direction is to make use of meshes aligned with the magnetic field lines \cite{dudson2015bout++, hoelzl2021jorek, jardin2010computational}. However, this can only apply for configurations in which the magnetic field lines form nested flux surfaces, which is not the case, for instance, for tokamak disruptions that can have volume-filling stochastic magnetic field lines \cite{freidberg1987ideal}. Mesh-alignment also places limitations on software and algorithm flexibility. Several other avenues have been explored focusing on increased accuracy of heat flux discretizations, thereby allowing for general, non-field-aligned meshes. Recently, authors have considered higher order discretizations \cite{green2022efficient} and adaptive mesh refinement \cite{vogl2022mesh}, where the mesh is refined specifically in areas where large error due to anisotropy is expected. One can also make use of an auxiliary variable, explicitly resolving the heat flux along magnetic field lines as a variable in its own right. This formulation has been shown to better resolve the heat flux after discretization, reducing the aforementioned parallel-to-perpendicular cross-contamination. One such auxiliary variable implementation, provided by Gunter et al. \cite{gunter2007finite}, uses a mixed finite element method, with the temperature field defined on a continuous Lagrange element of order $n$, and the auxiliary heat flux variable on a discontinuous Lagrange element of order $n-1$. Other approaches to reduce cross-pollution include asymptotic preserving methods \cite{degond2012asymptotic, deluzet2019two, jin1999efficient, narski2014asymptotic, yang2019preserving}, and first-order hyperbolic system methods \cite{chamarthi2019first} (and references therein). \ccrf{In the former, a scaled version of temperature fluctuations along field lines is considered as an auxiliary variable. In the latter, auxiliary variables are introduced which correspond to the temperature's (non-directional) gradient components together with pseudo-time advancement terms.}{i.}

For large sparse matrices that arise in finite element discretizations, almost all fast and parallelizable solvers are iterative and multilevel in nature. The basic idea is to combine a relatively cheap procedure to attenuate error on one ``level'' in the hierarchy (``relaxation''), with a coarse-representation that represents the complementary error in a smaller solution space (``coarse-grid correction''). Relaxation is typically a local procedure, such as Jacobi, Gauss-Seidel, or (overlapping) Schwarz, which attenuates high-frequency error in a geometric sense, and error associated with large eigenvalues in a spectral sense. The fundamental difficulty with developing fast solvers for strongly anisotropic diffusion is that there exists a large set of eigenmodes that are smooth in the parallel direction and high frequency in the perpendicular direction, which have small eigenvalues. Such modes cannot be effectively attenuated by standard local relaxation procedures due to the small associated eigenvalues, and are equally difficult to represent in a coarse solution space due to the high frequency character in the perpendicular direction.

One solution to the aforementioned problems is \emph{line-relaxation}, where one does a block relaxation aligned with one dimension. E.g., in 2D, for each discrete point in the $x$-dimension, one solves the fully coupled problem in the $y$-dimension. If the mesh is aligned with the anisotropy, one can perform line relaxation in the direction of anisotropy; for non-grid aligned anisotropy, one can use alternating line relax in grid-aligned directions, e.g., in $x$ and $y$ in 2D, or attempt to find an ordering of mesh points that is approximately aligned with the anisotropy. In 3D, things become more complicated, typically requiring \emph{plane} relaxation. Line relaxation was a key component of the multilevel solvers recently used in \cite{green2022efficient} for anisotropic heat flux in plasma. The downsides of line relaxation are its high cost, need for structured grids to define lines/planes on, preferably aligned with anisotropy, and in particular its poor performance in parallel. A number of papers have also tried to address anisotropic diffusion without line relaxation through algebraic solvers, primarily variations in AMG (e.g., see \cite{Brannick.2012,DAmbra.2013,Gee.2009,manteuffel2017root,Schroder.2012}). Although there has been some success in this regard, there is a disconnect between ``hard'' anisotropic problems considered in solvers papers, and accurate discretizations of the anisotropic heat flux that actually arises in, e.g., magnetic confinement fusion problems. Further, for general meshes and the implicit regime of interest for magnetic confinement fusion, there is also a lack of efficient solver strategies for the other aforementioned discretization approaches.

This paper conjointly addresses the challenges in discretization accuracy and solver efficiency for strongly anisotropic heat flux with open field lines. Following \cite{gunter2007finite}, we pose the anisotropic heat flux in mixed form using an auxiliary variable, which offers superior accuracy in the strongly anisotropic regime compared with direct discretizations of the second-order diffusion operator. This formulation also uncovers the purely anisotropic diffusion operator as a combination of two linear advection operators. The $2\times2$ block system resulting from a mixed discretization typically has mass matrices and the isotropic diffusion term on the diagonal, with the advection operators in the off-diagonal blocks. Although robust multilevel methods for (physically realistic) anisotropic diffusion remains a largely open question, significant progress has been recently made developing fast and robust AMG methods for linear advection based on an approximate ideal restriction (AIR) \cite{manteuffel2019nonsymmetric, manteuffel2018nonsymmetric}. Thus, rather than follow the traditional approach approximating a Schur complement that represents the anisotropic diffusion operator, we instead swap the rows of the block system so that the advection operators are on the diagonal, and apply a block triangular preconditioning, with the diagonal advection blocks inverted using AIR-AMG. Given the transport point of view of this solver scheme, it is well-constructed to deal with the anisotropy of the problem: the higher the anisotropy, the more efficient this approach is compared to classical ones. While this approach is similar to the aforementioned first-order hyperbolic system methods, the key difference is that here, we exploit the formulation of two advection operators for the purpose of solver efficiency.

The remainder of this work is structured as follows: in \Cref{sec_motivation} we introduce the problem and motivate our solver strategy. The proposed solver strategy has several implications on the discretization that must be taken into consideration. In particular, the advection blocks need to be square and invertible, which requires care in choice of finite element space, enforcement of boundary conditions, and constraints on the anisotropic field lines. In addition, we must choose a discretization that AIR is an effective solver for -- as for many problems, some FEM formulations are much easier to solve than others. \Cref{sec_spat_disc} describes a new DG-upwind-based spatial discretization of the auxiliary formulation that conforms with each of these constraints. The full solver strategy is then introduced in \Cref{sec:solvers}, including a spectral analysis of the (advection-)preconditioned Schur complement. Numerical results are presented in \Cref{sec_Numerical_results}. From an accuracy perspective, the new mixed DG discretization is shown to significantly outperform a non-mixed CG and DG formulation, as well as the mixed CG system developed in \cite{gunter2007finite}, and is \emph{fully applicable and demonstrated on closed field lines}. In addition, in the regime of large anisotropy and open field lines, the proposed solver strategy is shown to be fast in terms of iteration counts and wall-clock times, whereas other classical diffusion-based AMG schemes demonstrate very poor scaling or simply fail to converge. Finally, in Section \ref{sec_conclusion}, we review our results and discuss future work.
\pagebreak
\section{Motivation} \label{sec_motivation}
\subsection{Anisotropic Heat Transport} \label{subsec_model}

Following Gunter et al. ~\cite{gunter2007finite}, we formulate the evolution of temperature $T$ within a plasma through use of a heat flux $\mathbf{q}$, with forcing $S$:
\begin{subequations}\label{T_eqn_orig}
\begin{align}
\pp{T}{t} - \nabla \cdot \mathbf{q} & = S,  \\
\mathbf{q} &= \kappa_{\parallel} \nabla_{\parallel} T + \kappa_{\perp} \nabla_{\perp} T, \label{q_orig}
\end{align}
\end{subequations}
over some bounded domain $\Omega$, and where $\nabla_{\parallel}$ and $\nabla_{\perp}$ are the components of the gradient along and orthogonal to magnetic field lines, respectively:
\begin{equation}
\nabla_{\parallel}(\cdot) := \big(\vect{b}\cdot\nabla(\cdot)\big)\vect{b}, \;\;\;\;\;\;\;\; \nabla_\perp:= \nabla - \nabla_\parallel,
\end{equation}
for $\vect{b} = \vect{B}/|\vect{B}|$, where $\vect{B}$ is the magnetic field along which the anisotropy is defined. Similarly, $\kappa_{\parallel}$ denotes the conductivity parallel to magnetic field lines and $\kappa_{\parallel}$ denotes the conductivity perpendicular to field lines. They may depend on the magnetic field and the temperature and for implicit discretizations, the resulting nonlinear problem can be solved, for example, with a quasi-Newton iteration scheme, using the given iteration's known value for $T$ in the conductivities. For simplicity, in what follows, we assume constant conductivities and note that the framework follows equally for non-constant ones (up to the required nonlinear iterative procedure). \ccbf{It remains to specify boundary conditions. In magnetic confinement fusion applications, heat transport is often considered together with pure Dirichlet boundary conditions for $T$:
\begin{equation}
T|_{\partial \Omega} = T_\text{BC},
\end{equation}
and we will assume such boundary conditions in the main sections of this work. Alternatively, one can consider so-called sheath-boundary conditions (see e.g.~\cite{schneider2006plasma,stangeby2000plasma}), which are of Neumann type. In  \Cref{appendix:Neumann_BC}, we describe an adjustment for such boundary conditions, given the Dirichlet boundary condition-based space discretization to be introduced in Section \ref{sec_spat_disc}.}{ii.}

Heat flux in the context of magnetic confinement fusion is highly anisotropic, diffusing significantly faster along magnetic field lines as opposed to perpendicular to them. As mentioned in the introduction, we can expect to see $\frac{\kappa_{\parallel}}{\kappa_{\perp}} \sim 10^{10}$ and larger. \ccbf{Next to splitting the heat flux into a purely parallel and purely perpendicular component, it is also possible to split it into an isotropic part, and a remaining purely anisotropic part. In our case, the latter part is in the parallel direction, and its conductivity is given by}{2} $\kappa_\Delta \coloneqq \kappa_\parallel - \kappa_\perp$. Heat conduction in the plasma is isotropic when $\kappa_\Delta = 0$, and is predominately anisotropic if $\kappa_\Delta \gg \kappa_\perp$. Defining $\kappa_\Delta$ in this way allows us to represent the heat flux as
\begin{equation}
\mathbf{q}  = \kappa_{\Delta} \nabla_{\parallel} T + \kappa_{\perp} \nabla T.
\end{equation}
Considering the directional derivative $\vect{b} \cdot \nabla T$ as an auxiliary variable, we may rewrite \eqref{T_eqn_orig} as
\begingroup
\begin{subequations} \label{T_eqn_full}
\begin{align}
\pp{T}{t} - \sqrt{\kappa_\Delta} \nabla \cdot(\vect{b} \zeta) - \nabla \cdot (\kappa_\perp \nabla T) & = S, \label{T_eqn_full_T} \\
\zeta & = \sqrt{\kappa_\Delta} \vect{b} \cdot \nabla T, \label{T_eqn_full_zeta}
\end{align}
\end{subequations}
\endgroup
where we split up the factor of $\kappa_\Delta$ evenly between equations \eqref{T_eqn_full_T} and \eqref{T_eqn_full_zeta} in view of the solver strategy to be presented below. In this form, the PDE system has a more conventional advection-diffusion character. Gunter et al.~\cite{gunter2007finite} demonstrated that this mixed system -- solved for $T$ and $\zeta$ -- may lead to higher accuracy than a single-variable model which solves for $T$ only, by resolving the parallel heat flux term $\zeta$ more accurately.

Finally, we note the equation's two limiting cases in order to motivate the solver strategy presented in this work: for $\kappa_\parallel = \kappa_\perp$, we find that $\kappa_\Delta = 0$ and the system reduces to a (forced) heat equation given by
\begin{align}
\pp{T}{t} = \nabla \cdot (\kappa_\perp \nabla T) + S, \label{T_eqn_iso}
\end{align}
and we can expect classical multigrid methods to perform well when solving implicit discretizations of \eqref{T_eqn_iso}. However, for $\kappa_\perp = 0$, we instead have a purely anisotropic system of equations
\begingroup
\begin{subequations} \label{T_eqn_aniso}
\begin{align}
\pp{T}{t} - \sqrt{\kappa_\Delta} \nabla \cdot(\vect{b} \zeta) &= S, \\
\zeta &= \sqrt{\kappa_\Delta}\vect{b} \cdot \nabla T,
\end{align}
\end{subequations}
\endgroup
in which the dynamics are governed by two transport operators.

\subsection{Solver considerations} \label{sec_solver_considerations}

Written in block operator form, \eqref{T_eqn_full} can be formulated as
\begin{equation}\label{mat-cts}
\begin{pmatrix}
  \pp{}{t} + \kappa_\perp \Delta && -\sqrt{\kappa_\Delta} \nabla \cdot \big(\mathbf{b} (\cdot)\big) \\
  -\sqrt{\kappa_\Delta} \mathbf{b} \cdot \nabla && I
\end{pmatrix}
\begin{pmatrix}
  T \\
  \zeta
\end{pmatrix}
=
\begin{pmatrix}
  S \\
  0
\end{pmatrix},
\end{equation}
where $\Delta$ denotes the Laplacian, and $I$ the identity operator. Upon discretizing in space with a suitable finite element method and applying implicit integration in time (including diagonally implicit Runge Kutta or linear multistep methods), the resulting discrete linear systems may take the form
\begin{equation}\label{mat-discr}
\begin{pmatrix}
  \frac{1}{\Delta t}M_T + \kappa_\perp L && \sqrt{\kappa_\Delta} G_b^T \\
  -\sqrt{\kappa_\Delta}G_b && M_\zeta
\end{pmatrix}
\begin{pmatrix}
  T^{n+1} \\
  \zeta^{n+1}
\end{pmatrix}
=
\begin{pmatrix}
  F_T \\
  F_\zeta
\end{pmatrix},
\end{equation}
where $M_T$ and $M_\zeta$ are mass matrices, $L$ denotes the discrete Laplacian, $G_b$ the discrete scalar-form transport operator, and $F_T$, $F_\zeta$ correspond to forcing terms as well as terms on the known time level $n$.

Consider the limiting case of a steady-state (\ccrf{i.e.~}{1}no time derivative) and purely anisotropic problem ($\kappa_\perp = 0$). Plugging this value of $\kappa_\perp$ into \eqref{mat-discr} and swapping the block rows in the corresponding equations yields the equivalent block-triangular system
\begin{equation}\label{mat-discr_steady}
\begin{pmatrix}
  -\sqrt{\kappa_\Delta}G_b && M_\zeta \\
  0 && \sqrt{\kappa_\Delta} G_b^T
\end{pmatrix}
\begin{pmatrix}
  T^{n+1} \\
  \zeta^{n+1}
\end{pmatrix}
=
\begin{pmatrix}
  F_\zeta \\
  F_T
\end{pmatrix}.
\end{equation}
Here, we have decoupled the implicit purely anisotropic diffusion equation into the solution of two successive linear advection equations. Although it may seem unusual to take a parabolic problem and reformulate in terms of hyperbolic operators, in this case there is a fundamental lack of efficient, parallelizable iterative solvers for the parabolic purely anisotropic setting, while there do exist fast AMG methods for linear advection \cite{manteuffel2019nonsymmetric, manteuffel2018nonsymmetric}. Note that we opted for a multiplicatively even split of $\kappa_\Delta$ between $G_b$ and $G_b^T$ in view of future work considering temperature dependent conductivities. While the split is not necessary for the constant values of $\kappa_\Delta$ considered here, it may avoid that one solve is more challenging than the other when such a temperature dependency is included.

The above limiting case motivates the more general solver strategy proposed in this paper (and detailed in \Cref{sec:solvers}), focused on inverting transport operators rather than directly trying to solve the highly anisotropic diffusion equations. For this framework to be viable, we must construct the transport operator to be invertible, which can be achieved by (i) choosing the same finite element space for $T$ and $\zeta$, and (ii) assuming \emph{open} magnetic field lines throughout the domain, that is, the magnetic field is acyclic. The former ensures that $G_b$ is square, and the latter ensures that the operator is nonsingular -- in contrast, if the magnetic field has closed field lines, then constant values of $T$ along any closed field line constitute a non-trivial kernel. In addition, we must choose a transport discretization that is amenable to the fast reduction based AMG solver, AIR, developed in \cite{manteuffel2019nonsymmetric, manteuffel2018nonsymmetric}. Because AIR is reduction based, it is most effective on discretizations with low matrix connectivity, a property \ccrf{particularly}{2} enjoyed by upwind discontinuous Galerkin discretizations, where matrix connectivity only occurs to immediately adjacent elements that are upwind with respect to the magnetic field. The following section constructs a suitable discretization, and the solver strategy is expanded in \Cref{sec:solvers}.

\section{Finite element discretization}\label{sec_spat_disc}

To solve the anisotropic diffusion system \eqref{T_eqn_full}, we discretize the transport term
\begin{equation}
  \nabla \cdot \big(\sqrt{\kappa_\Delta} \vect{b} (\cdot)\big),
\end{equation}
using the classical DG-upwind method (see e.g.~\cite{kuzmin2010guide}). For this purpose, we start with an operator of the form
\begingroup
\addtolength{\jot}{2mm}
\begin{align}
  L_b (\theta ; \phi) := &- \left \langle \sqrt{\kappa_\Delta}\theta \vect{b}, \nabla \phi \right \rangle + \int_\Gamma \llbracket \sqrt{\kappa_\Delta} \phi \vect{b} \cdot \mathbf{n} \rrbracket \tilde{\theta} \;\text{d}S +\int_{\partial \Omega_{\text{out}}} \sqrt{\kappa_\Delta} \phi \vect{b} \cdot \mathbf{n} \theta \;\text{d} S, \label{transport_out}
\end{align}
\endgroup
for $L^2$-inner product $\langle \cdot, \cdot \rangle$, and any functions $\theta$, $\phi$ in the given $k^{th}$ polynomial order DG space $\mathbb{V}^{\text{DG}}_k(\Omega)$, which we will specify further below in this section. The set $\partial \Omega_{\text{out}}$ denotes the outflow subset of the boundary $\partial \Omega$ relative to $\mathbf{B}$. Further, $\Gamma$ denotes the set of all interior facets of the mesh, and we applied jump and upwind facet operations defined by
\begin{equation}
\llbracket \psi \rrbracket \coloneqq \psi^+ - \psi^-, \hspace{2cm} \tilde{\psi} \coloneqq
\begin{cases}
\psi^+ \;\; \text{if } \mathbf{b}^+ \cdot \mathbf{n}^+ < 0,\\
\psi^- \;\; \text{otherwise},
\end{cases}
\end{equation}
for any scalar field $\psi$. $\mathbf{n}$ denotes the facet normal vector, and the two sides of each mesh facet are arbitrarily denoted by $+$ and $-$ (and hence $\mathbf{n}^+ = - \mathbf{n}^-$). Given \eqref{transport_out}, we then define two transport operators according to
\begingroup
\addtolength{\jot}{2mm}
\begin{subequations} \label{transport_T_zeta}
\begin{align}
&L_{b, T}(\zeta_h; \phi) \coloneqq L_b(\zeta_h; \phi) + \int_{\partial \Omega_{\text{in}}} \sqrt{\kappa_\Delta} \phi (\vect{b} \cdot \mathbf{n}) \zeta_{\text{in}} \;\text{d} S & \forall \phi \in \mathbb{V}^{\text{DG}}_k, \label{transport_T}\\
&L_{b, \zeta}(\psi; T_h) \coloneqq L_b(\psi; T_h) - \int_{\partial \Omega_{\text{out}}} \sqrt{\kappa_\Delta} T_{\text{BC}} (\vect{b} \cdot \mathbf{n}) \psi \;\text{d} S & \forall \psi \in \mathbb{V}^{\text{DG}}_k, \label{transport_zeta}
\end{align}
\end{subequations}
\endgroup
where $\partial \Omega_{\text{in}}$ denotes the inflow part of the boundary, defined analogously to $\partial \Omega_{\text{out}}$, and $\zeta_{\text{in}}$ denotes a known value to be specified at the end of this section. Note that this setup of known inflow boundary integrals ensures that $L_{b,T}$ and $L_{b,\zeta}$ are non-singular for problems where all magnetic field lines are open. An intuition for this choice of operators and boundary integrals is given in \Cref{appendix:upwind_transport}.

Given these operators, we can formulate the system to be solved as finding $(T_h, \zeta_h) \in (\mathbb{V}^{\text{DG}}_k \times \mathbb{V}^{\text{DG}}_k)$ such that\!
\begingroup
\addtolength{\jot}{4mm}
\begin{subequations} \label{DG_upw_discr}
\begin{align}
&\left\langle \phi, \pp{T_h}{t} \right \rangle - L_{b, T}(\zeta_h; \phi) - \text{IP}(T_h; \phi) \nonumber\\
&\hspace{32mm} = -\! \int_{\partial \Omega}\! \kappa_{\text{BC}} \phi (T_h-T_{\text{BC}}) \;\text{d}S + \langle \phi, S \rangle &\text{\qquad $\forall \phi \in \mathbb{V}^{\text{DG}}_k$}, \label{DG_upw_discr_T} \\
&\left\langle \psi, \zeta_h \right \rangle + L_{b, \zeta}(\psi; T_h) =0 &\text{\qquad $\forall \psi \in \mathbb{V}^{\text{DG}}_k$},\label{DG_upw_discr_zeta}
\end{align}
\end{subequations}
\endgroup
for penalty parameter $\kappa_\text{BC}$, and where IP($T_h;\phi$) is the discretized Laplacian using an interior penalty method (see e.g.~\cite{burman2005unified}):
\begingroup
\addtolength{\jot}{4mm}
\begin{align}
\begin{split}
\text{IP}(T_h; \phi) = &- \left \langle \nabla \phi, \kappa_\perp \nabla T_h \right \rangle
+\int_\Gamma \left \llbracket T_h \right \rrbracket \left \lbrace \kappa_\perp \nabla \phi \right \rbrace \;\text{d}S
+ \int_\Gamma \left \llbracket \phi \right \rrbracket \left \lbrace \kappa_\perp \nabla T_h \right \rbrace \;\text{d}S\\
&\hspace{-12mm} \!+\!\! \int_{\partial \Omega}\!\! \kappa_\perp \mathbf{n} \cdot \nabla \phi(T_h-T_{\text{BC}}) \;\text{d}S
\!+\!\! \int_{\partial \Omega}\!\! \kappa_\perp \mathbf{n} \cdot \nabla T_h \phi \;\text{d}S \!-\!\! \int_\Gamma \frac{\kappa_\perp \kappa_p}{h_e} \left \llbracket \phi \right \rrbracket \left \llbracket T_h \right \rrbracket \;\text{d}S, \label{IP_term}
\end{split}
\end{align}
\endgroup
for any test function $\phi \in \mathbb{V}_k^\text{DG}$, penalty parameter $\kappa_p$, and where $h_e$ denotes the given facet's length scale. For the numerical results below,  we used as a length scale $\tfrac{|K^+| + |K^-|}{2|\partial K|}$, where $|K^{\pm}|$ denote the facet's adjacent cell's volumes, and $|\partial K|$ denotes the facet area. Further, the average operation is given by
\begin{equation}
\{\mathbf{u}\} \coloneqq \frac{1}{2} \left(\mathbf{u}^+ + \mathbf{u}^-\right)\cdot\mathbf{n}^+,
\end{equation}
for any vector field $\mathbf{u}$. Finally, the non-dimensional interior penalty parameter in $\text{IP}(T_h; \phi)$ is set to
\ccrf{
\begin{equation}
\kappa_p = 2,
\end{equation}
}{3}
in the numerical results section below. The discussion for the penalty parameter $\kappa_\text{BC}$ for exterior facets (first term on right-hand side of \eqref{DG_upw_discr_T}) is more intricate and we postpone it to the end of this section.

In the numerical results section below, we consider 3D extruded meshes based on prism cells. Such cells are of interest in magnetic confinement fusion applications, where tokamak meshes may be created starting from a triangular mesh for the poloidal plane (\ccrf{i.e.~}{1}a vertical cross section), which is then extruded to create a 3D mesh. While the extrusion would be curved in the context of tokamak meshes, here we restrict ourselves to straight extrusions for simplicity\footnote{In exploratory numerical tests, we found our framework to also work in curved meshes, albeit with an additional error that depends on the finite element space and mesh degrees in the direction of the curved extrusion.}. We then set $\mathbb{V}^{\text{DG}}_k(\Omega)=\mathbb{V}^{\text{DG}}_2(\Omega)$ to the corresponding second polynomial order discontinuous Galerkin space. For prisms, the latter space can be constructed via a finite element which is given by the tensor product of the standard DG finite element $dP_2$ for triangles, together with the analogous one for intervals. Further, we discretize $\mathbf{B}$ in the second order div-conforming Nedelec space $N{c_2}^f$. There is no specific reason to this choice of finite element space other than related work on magnetic confinement fusion discretizations by the authors. With regards to the choice of magnetic field space, we generally found the scheme's accuracy to be only significantly affected by the latter space's polynomial degree, which should be at least equal to the temperature space's polynomial degree. \ccrf{Alternatively, it is possible to set $\mathbf{B}$ as an expression to be evaluated at quadrature points; here we decided to set $\mathbf{B}$ as a finite element field since anisotropic diffusion is often considered within a discretization of some form of the magnetohydrodynamic equations, in which $\mathbf{B}$ is evolved discretely.}{4}
\begin{remark} \label{Remark:structure}
\textbf{Structural properties}. The non-discretized anisotropic diffusion operator in \ccrf{\eqref{T_eqn_full_T}}{5} acts to diffuse temperature and may create a heat flux across the domain's boundary in the direction of the magnetic field $\mathbf{B}$. This can easily be seen by taking the $L^2$-inner product of \ccrf{\eqref{T_eqn_full_T}}{5} with $T$, and setting $\kappa_\perp = 0$. After integrating by parts, this yields
\begin{equation}
  \frac{1}{2} \frac{\text{d}}{\text{d}t} \|T\|_2^2 = - \| \zeta \|_2^2 + \int_{\partial \Omega} \sqrt{\kappa_\Delta}\zeta \vect{b} \cdot \mathbf{n} T_{\text{BC}} \;\text{d}S, \label{L2_T_change_nondiscr}
\end{equation}
up to contributions due to the forcing $S$. Further, the anisotropic transport operator does not impact the system's total energy other than at the boundary. This can by seen by integrating equation \ccrf{\eqref{T_eqn_full_T}}{5} (with $\kappa_\perp = 0$) and applying integration by parts again. We then have
\begin{equation}
  \frac{\text{d}}{\text{d}t} \int_\Omega T  \;\text{d}x = \int_{\partial \Omega} \sqrt{\kappa_\Delta}\zeta \vect{b} \cdot \mathbf{n} \;\text{d}S, \label{L1_T_change_nondiscr}
\end{equation}
where in this context, the integral over $T$ plays the role of total energy.

Similarly, for the newly introduced DG discretization \eqref{DG_upw_discr} with the anisotropic part only (with $\kappa_\perp = 0$), we can set $\phi = T_h$ in \eqref{DG_upw_discr_T} and further $\psi = \zeta_h$ in \eqref{DG_upw_discr_zeta}, which leads to
\begingroup
\addtolength{\jot}{4mm}
\begin{align}
\begin{split}
  \frac{1}{2} \frac{\text{d}}{\text{d}t} \|T_h\|_2^2 = &- \| \zeta_h \|_2^2 - \int_{\partial \Omega} \kappa_\text{BC} T_h (T_h-T_{BC}) \;\text{d}S \\
  &+ \int_{\partial \Omega_{\text{in}}} \sqrt{\kappa_\Delta} \zeta_{\text{in}} (\vect{b} \cdot \mathbf{n}) T_h \;\text{d}S + \int_{\partial \Omega_{\text{out}}} \sqrt{\kappa_\Delta} \zeta_h (\mathbf{b} \cdot \mathbf{n}) T_{\text{BC}} \;\text{d}S,
  \label{discrete_diffusion}
\end{split}
\end{align}
\endgroup
where we substituted $L_b(\zeta_h, T_h)$ from $L_{b, T}(\zeta_h, T_h)$ in \eqref{DG_upw_discr_T} by the corresponding term from \eqref{DG_upw_discr_zeta}. This behavior is consistent in the sense that when $T_h$ is equal to the strong solution $T$, and $\zeta_\text{in}$ is set to the strong directional gradient $\sqrt{\kappa_\Delta} \mathbf{b} \cdot \nabla T$, these terms are equal to the right-hand side of \eqref{L2_T_change_nondiscr}. Note that the correct discrete diffusive behavior \eqref{discrete_diffusion} may initially seem counterintuitive, since we obtained it by combining two transport operators that apply upwinding in the same direction (along
$\mathbf{b}$). This can be resolved by observing that while we upwind the trial function for $L_{b, T}$ in
\eqref{transport_T}, we instead upwind the test function for $L_{b, \zeta}$ in \eqref{transport_zeta}. The latter in turn can be considered as downwinding the trial function. Altogether, the upwinded and downwinded trial functions then combine to the above diffusive behavior.

Finally, to obtain a discrete version of \eqref{L1_T_change_nondiscr}, we set $\phi = 1$ in \eqref{DG_upw_discr_T}, which leads to
\begingroup
\addtolength{\jot}{2mm}
\begin{align}
\frac{\text{d}}{\text{d}t} \int_\Omega T_h  \;\text{d}x =& \int_{\partial \Omega_{\text{out}}} \sqrt{\kappa_\Delta} \zeta_h(\vect{b} \cdot \mathbf{n}) \;\text{d} S
 + \int_{\partial \Omega_{\text{in}}} \sqrt{\kappa_\Delta} \zeta_{\text{in}} (\vect{b} \cdot \mathbf{n}) \;\text{d} S \label{energy_discr}\\
 &- \int_{\partial \Omega} \kappa_\text{BC} (T_h-T_{BC}) \;\text{d}S + \int_\Gamma \llbracket \sqrt{\kappa_\Delta} \mathbf{b} \cdot \mathbf{n}\rrbracket \tilde{\zeta}_h \;\text{d} S. \nonumber
\end{align}
\endgroup
For the boundary integrals, as before this is consistent with the non-discretized version \eqref{L1_T_change_nondiscr} when $T_h$ is equal to the strong solution $T$. However, this time, there is an additional interior facet integral term (last term in \eqref{energy_discr}), which need not vanish since $\mathbf{b} = \mathbf{B}/|\mathbf{B}|$ need not have continuous normal components, even if $\mathbf{B}$ is a div-conforming finite element field (\ccrf{i.e.~}{1}has continuous normal components across facets). However, as for the boundary integrals, the term is weakly consistent in the sense that for a strong solution with a continuous field $\mathbf{B}$, the term vanishes.
\end{remark}

\ccbf{At this point, we stress that the upwind formulation is not only key to the AIR-based solver approach, but also to the more accurate field development results to be presented below. The latter results are achieved through a better representation of the two directional derivative operators appearing in the anisotropic diffusion operator. Like the structural properties demonstrated in \Cref{Remark:structure}, this holds true independently of whether or not the given problem's magnetic field lines are open or closed. This work's restriction to open field lines is therefore a consequence of solver considerations, rather than accuracy ones.}{i.}

We end this section with a brief description of the time discretization to be used in the numerical results section, as well as terminology for the next section, the choice of $\zeta_\text{in}$, and a discussion on the penalty terms. The space discretization \eqref{DG_upw_discr} is coupled to an implicit midpoint rule. Given the fully discretized scheme, we denote the operator $L_b$ in matrix form as $-\sqrt{\kappa_\Delta} G_b$, where the change in sign and factorization of $\sqrt{\kappa_\Delta}$ are to aid the transport solver discussion below. Similarly, we denote the discrete Laplacian operator $\text{IP}$ in matrix form as $-\kappa_\perp L$. Further, we denote the $\mathbb{V}^{\text{DG}}_2$ mass matrix by $M$. Finally, for the implicit system, the boundary penalty term (first term on righ-hand side of \eqref{DG_upw_discr_T}) leads to a mass matrix defined along the boundary only, and we denote the latter by $M_{\text{BC}}$.

The need for specifying $\zeta_\text{in}$ as a value independent of the unknown $\zeta_h$ to be solved for is a shortcoming of the proposed space discretized formulation \eqref{DG_upw_discr}, which arises from employing DG-upwinded transport operators for both $T_h$ and $\zeta_h$. We therefore require known inflow boundary values for both fields for the operators to be invertible, while only values for either $T_h$ (Dirichlet) or $\zeta_h$ (Neumann) are given. One possible way to circumvent this is to set $\zeta_\text{in} = \sqrt{\kappa_\Delta} \mathbf{b} \cdot \nabla \bar{T}_h$, where $\bar{T}_h$ is the midpoint-in-time discretization of $T_h$. In particular, this is an expression independent of $\zeta_h$, and therefore does not change the resulting transport block $-\sqrt{\kappa_\Delta} G_b$ which we are required to invert. While we found this approach to work well in terms of solution accuracy, it adds a contribution with overall factor $\kappa_\Delta$ to the top-left block of the matrix system \eqref{mat-discr}, thereby rendering our transport operator based solver motivation \eqref{mat-discr_steady} invalid. Instead, we initially set the inflow flux $\zeta_\text{in}$ equal to the projection of $\sqrt{\kappa_\Delta} \mathbf{b} \cdot \nabla T|_{t=0}$ into $\mathbb{V}^{\text{DG}}_2$, and subsequently set it equal to the most recent solution for $\zeta_h$ of the resulting mixed system. In other words, $\zeta_{\text{in}}$ can be seen as an explicit time discretization of $\zeta_h$. Generally, we found this setup to be stable, likely due to the interior penalty term that weakly enforces the Dirichlet boundary conditions. We further found the resulting heat flux error at the inflow boundary generally to be small, including test cases where the heat flux varies in time at the latter boundary. In practice, this shortcoming can likely further be remedied at no additional cost by incorporating updates to the value of $\zeta_{\text{in}}$ in an outer nonlinear iterative procedure. As mentioned in \Cref{subsec_model}, the latter procedure is needed in magnetic confinement fusion simulations where the conductivities $\kappa_\perp$, $\kappa_\parallel$ depend on the temperature $T_h$.

Given this choice of $\zeta_\text{in}$, both additional boundary integrals in \eqref{transport_T_zeta} are absorbed into the right-hand side vectors $F_T$ and $F_\zeta$, respectively, leading to an overall discrete $2\times 2$ block system of the form
\begin{equation}\label{mat-discr-bdry-orig}
\begin{pmatrix}
  \frac{1}{\Delta t}M + \kappa_\perp L + M_\text{BC} && \sqrt{\kappa_\Delta} G_b^T \\
  -\sqrt{\kappa_\Delta}G_b && M
\end{pmatrix}
\begin{pmatrix}
  T_h^{n+1} \\
  \zeta_h^{n+1}
\end{pmatrix}
=
\begin{pmatrix}
  F_T \\
  F_\zeta
\end{pmatrix},
\end{equation}
whose structure is equal to \eqref{mat-discr} described in our motivation section \Cref{sec_solver_considerations}, up to the additional boundary penalty term. The latter term can be considered as part of the discrete Laplacian $L$, in which case the correct penalty parameter formulation is
\begin{equation}
\kappa_\text{BC} = \tilde{\kappa}_\text{BC} \frac{\kappa_\perp}{h_e}, \label{BC_penalty_L}
\end{equation}
for some non-dimensional parameter $\tilde{\kappa}_\text{BC}$. The top-left block of \eqref{mat-discr-bdry-orig} can then be grouped according to $\frac{1}{\Delta t}M + \kappa_\perp (L + \tilde{\kappa}_\text{BC} M_{\text{BC}, h_e^{-1}})$, where we factored $\tilde{\kappa}_\text{BC}\kappa_\perp$ out of $M_\text{BC}$ to highlight the matrix grouping. Further, we denote the remaining operator as $M_{\text{BC}, h_e^{-1}}$, noting that it contains a factor of $h_e^{-1}$. However, this may be impractical for the anisotropic limit in which $\kappa_\perp$ may be negligible, thereby not ensuring a strong enough weak enforcement of the boundary conditions. In particular, the transport operators $L_{b, T}$ and $L_{b, \zeta}$ alone are not sufficient to enforce the boundary conditions. We therefore instead consider $M_\text{BC}$ independently of $L$, and use a penalty parameter of the form
\begin{equation}
\kappa_\text{BC} = \tilde{\kappa}_\text{BC} \frac{h_e}{\Delta t}, \label{BC_penalty_M}
\end{equation}
where again $\tilde{\kappa}_\text{BC}$ is non-dimensional. In the numerical results section below, we set $\tilde{\kappa}_\text{BC} = 20$. Given the factor of $\tfrac{1}{\Delta t}$, we group $M_\text{BC}$ together with the mass matrix, leading to an overall mixed system
\begin{equation}\label{mat-discr-bdry}
\begin{pmatrix}
  \frac{1}{\Delta t}(M + \tilde{\kappa}_\text{BC} M_{\text{BC}, h_e}) + \kappa_\perp L && \sqrt{\kappa_\Delta} G_b^T \\
  -\sqrt{\kappa_\Delta}G_b && M
\end{pmatrix}
\begin{pmatrix}
  T_h^{n+1} \\
  \zeta_h^{n+1}
\end{pmatrix}
=
\begin{pmatrix}
  F_T \\
  F_\zeta
\end{pmatrix},
\end{equation}
where again we factored $\tilde{\kappa}_\text{BC}/\Delta t$ out of $M_\text{BC}$ to highlight the matrix grouping. This time, we denote the remaining matrix as $M_{\text{BC}, h_e}$, since it contains a factor of $h_e$.

\section{Transport solvers for anisotropic diffusion}\label{sec:solvers}

For standard block diagonal, triangular, and LDU preconditioners applied to block $2\times 2$ linear systems, convergence of Krylov and fixed-point iterations are fully defined by the preconditioning of the Schur complement \cite{southworth2020fixed}. In the standard form \eqref{mat-discr} of an auxiliary-variable discretization of \eqref{T_eqn_orig} or \eqref{T_eqn_full}, the natural (1,1) Schur complement is effectively a discrete representation of the continuous anisotropic diffusion equation; as discussed in \Cref{sec:intro}, it is very difficult to construct effective preconditioners for such equations. Building on the discussion in \Cref{sec_solver_considerations}, we reorder the discrete $2\times 2$ block system in \eqref{mat-discr-bdry} to take the equivalent form
\begin{equation}\label{mat-discr_reorder}
\begin{pmatrix}
  -\sqrt{\kappa_\Delta}G_b && M\\
  \frac{1}{\Delta t}(M  + \tilde{\kappa}_\text{BC} M_{\text{BC}, h_e}) + \kappa_\perp L&& \sqrt{\kappa_\Delta} G_b^T
\end{pmatrix}
\begin{pmatrix}
  T^{n+1} \\
  \zeta^{n+1}
\end{pmatrix}
=
\begin{pmatrix}
  F_\zeta \\
  F_T
\end{pmatrix}.
\end{equation}
In doing this reordering, the (2,2) Schur complement defining convergence takes the form (note, the (1,1) Schur complement is the negative transpose of the (2,2))
\begingroup
\addtolength{\jot}{2mm}
\begin{align}
  S_{22} \coloneqq&  \sqrt{\kappa_\Delta}  G_b^T +
    \big(\tfrac{1}{\Delta t}(M  + \tilde{\kappa}_\text{BC} M_{\text{BC}, h_e})\big)
    (\sqrt{\kappa_\Delta}G_b)^{-1}M + \kappa_\perp L(\sqrt{\kappa_\Delta}G_b)^{-1}M. \label{eq:schur}
\end{align}
\endgroup
Note that $M_{\text{BC}, h_e}$ can be seen as a mass matrix over boundary facets, which is scaled by a factor of $h_e$. Altogether, we expect its eigenspectrum to be akin to that of a mass matrix over cells. For simplicity, we further assume $\tilde{\kappa}_{BC}$ to be $\mathcal{O}(1)$, so that from the point of view of the spectral analysis to follow, $M  + \tilde{\kappa}_\text{BC} M_{\text{BC}, h_e}$ can be replaced by $M$ without a qualitative change in the analysis. We then have
\begin{equation}\label{eq:schur_simplified}
  S_{22} \approx \sqrt{\kappa_\Delta} G_b^T +
    \big(\tfrac{1}{\Delta t}M + \kappa_\perp L\big)
    (\sqrt{\kappa_\Delta}G_b)^{-1}M.
\end{equation}

This section provides a spectral analysis on preconditioning \eqref{eq:schur_simplified} by inverting the diagonal transport block to support the proposed strategy. The analysis is based on the observation that the preconditioned Schur complement is similar via a simple mass-matrix similarity transform to the identity plus two symmetric positive definite (SPD) operators. By nature of this similarity, we expect spectral analysis of the resulting SPD operators to provide excellent measures of preconditioning (a property which is not always the case with advective operators and inverses). The purely anisotropic case is discussed in \Cref{sec:solvers:aniso} and the general case in \Cref{sec:solvers:gen}. We then use a Lanczos algorithm to compute eigenvalue bounds for the preconditioned Schur complement for realistic problems in \Cref{sec:solvers:eigs} and demonstrate consistency with the theory. Finally, we detail the solver in practice in \Cref{sec:solver:solver}.

\subsection{The purely anisotropic case}\label{sec:solvers:aniso}

Recall in this paper we are particularly interested in the strongly anisotropic regime, where $\kappa_\Delta \gg \kappa_\perp$. If we consider the limit of $\kappa_\perp = 0$, the approximate Schur complement \eqref{eq:schur_simplified} simplifies to
\begin{equation}
  S_{22} \coloneqq \sqrt{\kappa_\Delta} G_b^T +
    \tfrac{1}{\Delta t}M (\sqrt{\kappa_\Delta}G_b)^{-1}M.
\end{equation}
Now consider preconditioning $S_{22}$ with the diagonal transport block
\begin{equation}\label{eq:prec-S}
  (\sqrt{\kappa_\Delta} G_b^T)^{-1}S_{22} = \
    I +
    \tfrac{1}{\Delta t\kappa_\Delta}G_b^{-T}M
    G_b^{-1}M.
\end{equation}
Observe that our preconditioned Schur complement is a perturbation of the identity, and (formally) similar to the SPD operator
\begin{equation}\label{eq:pert}
  M^{1/2}(\sqrt{\kappa_\Delta} G_b^T)^{-1}S_{22}M^{-1/2} =
  I + \tfrac{1}{\Delta t\kappa_\Delta}M^{1/2}G_b^{-T}M G_b^{-1}M^{1/2}.
\end{equation}
Next, we uncover the eigenvalue structure of the second term on the right-hand side of \eqref{eq:prec-S}.
Observe that $M^{-1} G_b^T M^{-1} G_b$ corresponds to a discrete anisotropic diffusion operator, which is similar to $M^{-1} G_b M^{-1} G_b^T$. Moreover, recall that $\mathbf{b}$ is a unit vector, and the eigenvalues of a discrete diffusion operator span a range of $[1/c_2,1/c_1h^2]$, for mesh spacing $h$ and $\mathcal{O}(1)$ constants $c_1,c_2$, independent of $h$. Thus, the inverse $G_b^{-T}M G_b^{-1}M$ that appears in \eqref{eq:prec-S}, and by similarity the second term of \eqref{eq:pert}, is an SPD operator with eigenvalues in the range $[c_1h^2, c_2]$. Altogether, this implies that the preconditioned Schur complement \eqref{eq:prec-S} has real-valued eigenvalues that are bounded below by one and lie in the range
\begin{equation}\label{eq:eig_bound1}
  \sigma \left[(\sqrt{\kappa_\Delta} G_b^T)^{-1}S_{22}\right]
    \subset 1 + \tfrac{1}{\Delta t\kappa_\Delta}[c_1 h^2, c_2].
\end{equation}
In particular, for reasonable values of $\Delta t$ and particularly for large $\kappa_\Delta$, we achieve excellent preconditioning of the Schur complement by simply inverting the diagonal advection blocks, largely independent of any additional constants scalings that depend on discretization. Of course, in the limit of a steady-state problem, $\Delta t\to\infty$, we have exactness, as used to motivate the solver strategy in \Cref{sec_solver_considerations}.

\subsection{The general case}\label{sec:solvers:gen}

For $\kappa_\perp > 0$, the analysis above largely holds, but we have additional perturbations in our preconditioned Schur complement, where
\begingroup
\addtolength{\jot}{2mm}
\begin{align}\label{eq:prec-S2:no_pert}
  (\sqrt{\kappa_\Delta} G_b^T)^{-1}S_{22}
  = I + \tfrac{1}{\Delta t\kappa_\Delta}G_b^{-T}&M G_b^{-1}M + \tfrac{\kappa_\perp}{\kappa_\Delta}G_b^{-T}LG_b^{-1}M.
\end{align}
\endgroup
As before, the above operator is similar to an SPD operator, this time given by
\begingroup
\addtolength{\jot}{2mm}
\begin{align}\label{eq:prec-S2}
\begin{split}
  M^{1/2}(\sqrt{\kappa_\Delta} G_b^T)^{-1}S_{22}M^{-1/2}
  = I + \tfrac{1}{\Delta t\kappa_\Delta}M^{1/2}G_b^{-T}&M G_b^{-1}M^{1/2} \\
  + \tfrac{\kappa_\perp}{\kappa_\Delta}M^{1/2}G_b^{-T}&LG_b^{-1}M^{1/2}.
\end{split}
\end{align}
\endgroup
Note that both of the perturbation terms in \eqref{eq:prec-S2} are SPD, and thus the preconditioned Schur complement again has positive real-valued eigenvalues, bounded below by one. As done already for the second-to-last term in \eqref{eq:prec-S2:no_pert}, we aim to establish a connection of the last term with a diffusion operator. Observe that by similarity, we have
\begin{equation}
G_b^{-T}LG_b^{-1}M \sim G_b^{-1}M G_b^{-T}L = \big(G_b^{-1}M G_b^{-T}M\big) \big(M^{-1}L\big), \label{eq:anti_diff_diff}
\end{equation}
where the latter term, $M^{-1}L$, corresponds to a discrete representation of the isotropic diffusion operator, and the first term, $G_b^{-1}M G_b^{-T}M$, corresponds to the inverse of the anisotropic diffusion operator as discussed in the previous section. To understand the product of these operators, suppose we are in 2D and $\mathbf{b}$ is grid-aligned, say in the $x$-dimension. Then \eqref{eq:anti_diff_diff} is a discrete representation of $\partial_{xx}^{-1}(\partial_{xx} + \partial_{yy}) = 1 + \partial_{xx}^{-1}\partial_{yy}$. On a unit domain with zero Dirichlet boundary conditions, the eigenfunctions and eigenvalues of the Laplacian are $\{u_{jk}, \pi^2(j^2+k^2)\}$, where $u_{jk} = 2\sin(j\pi x)\sin(k\pi y)$, and similarly for $\partial_{xx}$ only in the $(j,x)$ dimensions. Thus, if we consider the highest frequency mode represented by our grid in the $y$-dimension, $k = 1/h$, and the smoothest mode in the $x$-dimension, $j = 1$, we have the eigenpair
\begin{equation}\label{eq:fourier}
  (1 + \partial_{xx}^{-1}\partial_{yy})2\sin(\pi x)\sin(\tfrac{1}{h}\pi y)
    = \left(1 + \tfrac{1}{h^2}\right)2\sin(\pi x)\sin(\tfrac{1}{h}\pi y).
\end{equation}
A complete analysis generalized to specific FEM discretizations, non-unit-square domains, more complex boundary conditions, or non-grid aligned magnetic field $\mathbf{b}$ is a more complex endeavor and outside the scope of this work. However, a similar result will still hold, namely that modes that are highly oscillatory orthogonal to $\mathbf{b}$ and smooth in the direction of $\mathbf{b}$ will typically scale like $1+c_3/h^2$, with the $c_3/h^2$ arising due to oscillatory functions that are not damped by the anisotropic inverse. Altogether, we then have an estimate for the spectrum contained in the real-valued subset
\begin{equation}\label{eq:gen-eig-bound}
  \sigma\left(\tfrac{\kappa_\perp}{\kappa_\Delta}M^{1/2}G_b^{-T}LG_b^{-1}M^{1/2}\right)
    \subset \left(0, \frac{\kappa_\perp c_3}{\kappa_\Delta h^2}\right],
\end{equation}
for some constant $c_3$.

To conclude, the analysis in the previous two subsections indicates that for sufficiently small ratios of
\begin{equation}
\frac{c_1}{\Delta t \kappa_\Delta} \lessapprox 1,  \hspace{15mm} \frac{\kappa_\perp c_3}{\kappa_\Delta h^2} \lessapprox 1, \label{parameter_regimes}
\end{equation}
for constants $c_1$, $c_3$ independent of $\Delta t, \kappa_\Delta, h$, we expect excellent preconditioning of the Schur complement. Numerical results in \Cref{sec_Numerical_results} confirm these observations in practice, with rapid convergence obtained in solving \eqref{mat-discr_reorder}. On the other hand, for more isotropic regimes $\mathcal{O}(\kappa_\parallel) \approx \mathcal{O}(\kappa_\perp)$, and for setups with smaller time steps or $\kappa_\Delta$ such that $\mathcal{O}(\Delta t) \approx \mathcal{O}(\kappa_\Delta^{-1})$, we expect our solver performance to deteriorate since the transport operator alone will no longer be a good preconditioner for the Schur complement. Future work will consider more general methods that are also robust in the isotropic or mixed regimes. However, we note that in strongly magnetized plasmas, we are primarily concerned with extreme anisotropies for which the methods proposed here are robust.

Finally, in terms of the time step and therefore the first ratio in \eqref{parameter_regimes}, we set $\Delta t = 10^{-3}$ for the numerical tests in \Cref{sec_Numerical_results}. In comparison, the local diffusive time scale $h^2/\kappa_\parallel$ -- which is an important quantity in tokamak applications -- lies in an approximate range $h^2/\kappa_\parallel \in [10^{-4}, 10^{-13}]$ for the resolutions and parallel conductivities considered below. However, the latter scale is not the only factor determining the choice of time step and \ccrf{by}{6} extension the scheme's temporal accuracy. Other factors include the problem's nonlinearity in the case of $T$-dependent conductivities and the temperature profile. Overall, we found our fixed choice of time step to work well in terms of accuracy for the test cases considered below.

\subsection{Eigenvalue bounds in practice}\label{sec:solvers:eigs}

To further investigate the parameter regimes \eqref{parameter_regimes}, we include a brief numerical study to analyze the largest eigenvalue of the components of the preconditioned Schur complement $(\sqrt{\kappa_\Delta}G_b^T)^{-1}S_{22}$, with $S_{22}$ given by \eqref{eq:schur}. In the following, we denote the second-to-last term in \eqref{eq:schur} by $S_{22, M}$, and the last one by $S_{22, L}$. We have tested each operator using the mesh described in \Cref{appendix:base_mesh}, as well as a periodically extruded 3D version thereof. We consider 3 levels of mesh refinement, with 2 extruded cells in 3D for the lowest refinement. In each refinement, the resolution is doubled in each dimension. For each resolution, we tested the behavior for differing choices of time step $\Delta t \in [10^{-3}, 10^{-2}]$ and $\frac{\kappa_{\parallel} }{\kappa_{\perp}} \in [10^3, 10^6, 10^9]$, with $\kappa_\perp = 1$. Finally, the magnetic field $\mathbf{B}$ is defined as the curved and slanted one used in \Cref{sec_efficiency} below. The eigenvalues are computed by implementing the preconditioned Schur complement's components as actions applied to an input vector $\mathbf{x}_\text{in}$. The actions are then passed to an implicitly restarted Lanczos method \cite{lehoucq1998arpack} as implemented in SciPy\footnote{The number of Lanczos vectors is capped at 10 for the highest resolution 3D, together with a relative tolerance of 10\% for computational speedup, noting that the latter resolution leads to approximately $4.5\times10^5$ degrees of freedom.}. The resulting largest eigenvalues for the 2D and 3D setups are given in Tables \ref{Eigen_table} and \ref{Eigen_table_3D}, respectively; a larger leading eigenvalue will lead to a less efficient solver.
\begin{table}
\begin{center}
{\setlength{\extrarowheight}{3pt}
\setlength{\tabcolsep}{3.7pt}
\begin{tabular}{|c|c||c|c|c||c|c|c||c|c|c|}
 \hline
$\Delta t$&${\kappa_\parallel}/{\kappa_\perp}$ & \multicolumn{3}{c||}{$(\sqrt{\kappa_\Delta} G_b^T)^{-1}S_{22, M}$} & \multicolumn{3}{c||}{$(\sqrt{\kappa_\Delta} G_b^T)^{-1}S_{22, L}$} &  \multicolumn{3}{c|}{$(\sqrt{\kappa_\Delta} G_b^T)^{-1}S_{22} - I$}\\
 \hline
& $10^3$ & 7.8e0 & 4.6e0 & 2.9e0 & 5.6e-1 & 1.4e0 & 3.7e0 & 7.8e0 & 4.7e0 & 4.3e0\\
$10^{-3}$ & $10^6$ & 7.8e-3 & 4.6e-3 & 2.9e-3 & 5.6e-4 & 1.4e-3 & 3.7e-3 & 7.8e-3 & 4.7e-3 & 4.3e-3\\
& $10^9$ & 7.8e-6 & 4.6e-6 & 2.9e-6 & 5.6e-7 & 1.4e-6 & 3.7e-6 & 7.8e-6 & 4.7e-6 & 4.3e-6\\
  \hline
& $10^3$ & 7.8e-1 & 4.6e-1 & 2.9e-1 & 5.6e-1 & 1.4e0 & 3.7e0 & 8.0e-1 & 1.5e0 & 3.7e0\\
$10^{-2}$& $10^6$ & 7.8e-4 & 4.6e-4 & 2.9e-4 & 5.6e-4 & 1.4e-3 & 3.7e-3 & 8.0e-4 & 1.5e-3 & 3.7e-3\\
& $10^9$ & 7.8e-7 & 4.6e-7 & 2.9e-7 & 5.6e-7 & 1.4e-6 & 3.7e-6 & 8.0e-7 & 1.5e-6 & 3.7e-6\\
 \hline
\end{tabular}}
\caption{Largest eigenvalues of preconditioned Schur complement components for 2D mesh. $S_{22, M}$ and $S_{22, L}$ correspond to second and third terms on right-hand side of \eqref{eq:schur}, respectively. The three values for each given anisotropy ratio correspond to the three spatial refinement levels.} \label{Eigen_table}
\end{center}
\end{table}
\begin{table}
\begin{center}
{\setlength{\extrarowheight}{3pt}
\setlength{\tabcolsep}{3.7pt}
\begin{tabular}{|c|c||c|c|c||c|c|c||c|c|c|}

 \hline
$\Delta t$&${\kappa_\parallel}/{\kappa_\perp}$ & \multicolumn{3}{c||}{$(\sqrt{\kappa_\Delta} G_b^T)^{-1}S_{22, M}$} & \multicolumn{3}{c||}{$(\sqrt{\kappa_\Delta} G_b^T)^{-1}S_{22, L}$} &  \multicolumn{3}{c|}{$(\sqrt{\kappa_\Delta} G_b^T)^{-1}S_{22} - I$}\\
 \hline
                 & $10^3$ & 6.6e2 & 3.9e2 & 2.4e2 & 3.6e1 & 1.0e2 & 2.7e2 & 6.6e2 & 4.0e2 & 3.4e2\\
$10^{-3}$ & $10^6$& 6.6e-1 & 3.9e-1 & 2.4e-1 & 3.6e-2 & 1.0e-1 & 2.7e-1 & 6.6e-1 & 4.0e-1 & 3.4e-1\\
                & $10^9$ & 6.6e-4 & 3.9e-4 & 2.4e-4 & 3.6e-5 & 1.0e-4 & 2.7e-4 & 6.6e4 & 4.0e-2 & 3.4e-4\\
  \hline
                & $10^3$ & 6.6e1 & 3.9e1 & 2.4e1 & 3.6e1 & 1.0e2 & 2.7e2 & 6.6e1 & 1.1e2 & 2.8e2 \\
$10^{-2}$ & $10^6$& 6.6e-2 & 3.9e-2 & 2.4e-2 & 3.6e-2 & 1.0e-1 & 2.7e-1 & 6.6e-2 & 1.1e-1 & 2.8e-1 \\
                & $10^9$ & 6.6e-5 & 3.9e-5 & 2.4e-5 & 3.6e-5 & 1.0e-4 & 2.7e-4 & 6.6e-5 & 1.1e-4 & 2.8e-4 \\
 \hline
\end{tabular}}
\caption{Largest eigenvalues of preconditioned Schur complement components for 3D mesh. $S_{22, M}$ and $S_{22, L}$ correspond to second and third terms on right-hand side of \eqref{eq:schur}, respectively. The three values for each given anisotropy ratio correspond to the three spatial refinement levels.} \label{Eigen_table_3D}
\end{center}
\end{table}

For fixed resolutions, the largest eigenvalues $\lambda^{\max}_M$ of $(\sqrt{\kappa_\Delta} G_b^T)^{-1}S_{22, M}$ and $\lambda^{\max}_{L}$ of $(\sqrt{\kappa_\Delta} G_b^T)^{-1}S_{22, L}$ evolve according to their leading factors in \eqref{parameter_regimes}\ccrf{, where only the former depends on $\Delta t$, $\lambda^{\max}[(\sqrt{\kappa_\Delta} G_b^T)^{-1}S_{22, M}] \sim (\kappa_\Delta \Delta t)^{-1}$. In particular, while we only considered two choices of time step, this allows for extrapolation to other choices of time steps. For instance, for $\Delta t = 10^{-6}$ and $\tfrac{\kappa_\parallel}{\kappa_\perp} = 10^6$, we expect a largest eigenvalue akin to the one for $\Delta t = 10^{-3}$ and $\tfrac{\kappa_\parallel}{\kappa_\perp} = 10^3$. The corresponding value for the 3D setup is $\lambda^{\max}_M > 1\text{e}2$ for our choices of spatial resolution, and we therefore do not expect our solver to perform well in such a parameter regime. In contrast, at $\tfrac{\kappa_\parallel}{\kappa_\perp} = 10^9$, we still expect $\lambda^{\max}_M \sim 10^{-3}$ and thereby good iterative convergence.}{7}

The evolution according to the leading factors also holds true for varying resolutions, except increasing the resolution additionally decreases $\lambda^{\max}_M$. Upon further investigation, we found this to be related to the inclusion of the boundary penalty term $M_\text{BC}$ in the Schur complement $S_{22}$. Without it, the value $\lambda^{\max}_M$ is independent of the spatial resolution as expected -- for instance, in 2D and for $\Delta t = 10^{-3}$, $\kappa_\parallel/\kappa_\perp = 10^3$, we get $\lambda^{\max}_M = 0.84$ for all three resolutions. Finally, we found that as we increase the spatial resolution, the largest eigenvalues that are obtained when including $M_\text{BC}$ converge to the one when the latter operator is excluded.

We also observe that either $\lambda^{\max}_M$ or  $\lambda^{\max}_{L}$ may be larger, depending on the choice of $h$ and $\Delta t$ (for fixed conductivities). Additionally, the observed preconditioned Schur complement's largest eigenvalues (ignoring the additive factor of $1$ from the identity) are slightly less than the sum $\lambda^{\max}_M + \lambda^{\max}_L$. Finally, we find both $(\sqrt{\kappa_\Delta} G_b^T)^{-1}S_{22, L}$ and $(\sqrt{\kappa_\Delta} G_b^T)^{-1}S_{22,M}$ to have considerably larger largest eigenvalues in 3D than in 2D. In further studies including different discretization and initial condition setups, we found this likely to be due to the degree of the magnetic field's alignment with respect to the 3D mesh, rather than discretization factors such as the choice of finite element
, type of cell and cell aspect ratio.

\subsection{The solver in practice}\label{sec:solver:solver}

In practice, we apply a block lower triangular preconditioner of the form
\begin{equation}
\begin{pmatrix}
  -\sqrt{\kappa_\Delta}G_b && \mathbf{0}\\
  \frac{1}{\Delta t}(M + \tilde{\kappa}_\text{BC} M_\text{BC}) + \kappa_\perp L&& \sqrt{\kappa_\Delta} G_b^T
\end{pmatrix}^{-1},
\end{equation}
as a preconditioner to an outer flexible GMRES Krylov iteration solving the system in \eqref{mat-discr_reorder}, with an outer tolerance set to $10^{-8}$. We use lower triangular as opposed to upper triangular so as to directly include the isotropic operator in our preconditioning, even if is it not being inverted (and experimental results have supported this, with lower triangular preconditioning converging 10-20\% faster than upper triangular). Following \cite{southworth2020fixed}, we do not use block LDU preconditioning, noting that it typically provides minimal improvement in convergence compared with block triangular, at the expense of an additional (approximate) solve of a diagonal block. Inner solves of
\begingroup
\begin{subequations} \label{inner_solve_AIR}
\begin{align}
&(\sqrt{\kappa_\Delta} G_b)^{-1}, \label{inner_solve_AIR_inv}\\
&(\sqrt{\kappa_\Delta} G_b)^{-T}, \label{inner_solve_AIR_inv_T}
\end{align}
\end{subequations}
are approximated using AIR-AMG \cite{manteuffel2019nonsymmetric,manteuffel2018nonsymmetric} as a preconditioner for right-preconditioned GMRES. We choose an inner relative residual stopping tolerance of $\min\{ 10^{-3}, 10^{-3}/\|\bs{b}\|\}$, where $\bs{b}$ is the current right-hand side. This is effectively enforcing both an absolute and relative residual tolerance. We do this because when using block preconditioning with inexact Schur complements, it is almost always more computationally efficient to use approximate rather than exact inner solves (e.g., Navier Stokes results in \cite{southworth2020fixed} suggest a tolerance of about $10^{-3}$). However, when $(\Delta t \kappa_\Delta)^{-1}$ or $\kappa_\perp/\kappa_\Delta$ is not trivially small (in our experience $\gtrapprox 10^{-4}$), the right-hand side provided to the lower triangular solve can be very large, e.g., $\mathcal{O}(10^4)$. In such cases, we find that a relative stopping tolerance of $10^{-3}$ doesn't even lead to a residual $<1$ in norm, and the outer iteration fails to converge. Thus, for each inner solver we enforce achieving at least three digits in both absolute and relative residual tolerance.

Block AIR is applied using the PyAMG library \cite{bell2023pyamg}, with the discretization matrices grouped by DG element blocks. Ruge-Stuben coarsening \cite{ruge1987algebraic} with second pass is applied on a block level using a classical strength of connection with tolerance $\theta_C = 0.01$ based on a hard minimum, $\{-a_{ij} \geq \theta_C \max_k |a_{ik}|\}$, and calculated on a block basis as implemented in PyAMG. No pre-relaxation is used and FFC-block-Jacobi post relaxation is applied with no damping. We use one-point interpolation and distance-one block AIR \cite{manteuffel2019nonsymmetric} with strong connections determined through a classical strength of connection as above with tolerance $\theta_R=0.25$.

\section{Numerical results} \label{sec_Numerical_results}
Having described the scheme and solver strategy, we next test their properties in terms of convergence and efficiency. To set the scene, we first describe other discretizations and a standard solver strategy to compare against, as well as implementation details.

\subsection{Comparison to other discretizations and solvers} \label{comp_disc_solver}
For the test cases below, we will consider three additional discretizations. This includes the \textbf{CG-primal system} (\ccrf{i.e.~}{1}without the auxiliary heat flux term) for $T$ in the standard second polynomial order continuous Galerkin space $\mathbb{V}^\text{CG}_2$. For prisms, the latter is constructed analogously to the DG case, using the tensor product of the standard CG element $P_2$ for triangles and the analogous 1D element. In this case, all gradient operations are well-defined, and the discretization is given by
  \begin{align} \label{primal_CG}
  \left\langle \eta, \pp{T_h}{t} \right\rangle \!+\! \langle \mathbf{b} \cdot \nabla \eta, \kappa_\Delta \mathbf{b} \cdot \nabla T_h \rangle \!+\! \langle \nabla \eta, \kappa_\perp, \nabla T_h \rangle = \langle \eta, S\rangle && \forall \eta \in \mathbb{V}^{\text{CG}}_2.
  \end{align}
Additionally, we consider the $\mathbb{V}^{\text{CG}}_2$-$\mathbb{V}^{\text{DG}}_1$ \textbf{mixed finite element formulation} based on Gunter et al. \cite{gunter2007finite}, which reads
  \begingroup
  \addtolength{\jot}{2mm}
  \begin{subequations} \label{mixed_CG}
  \begin{align}
  &\left\langle \eta, \pp{T_h}{t} \right\rangle + \langle \mathbf{b} \cdot \nabla \eta, \sqrt{\kappa_\Delta} \zeta_h \rangle + \langle \nabla \eta, \kappa_\perp \nabla T_h \rangle = \langle \eta, S \rangle & \forall \eta \in \mathbb{V}^{\text{CG}}_2,\\
  &\langle \phi, \zeta_h - \sqrt{\kappa_\Delta} \mathbf{b} \cdot \nabla T_h \rangle = 0 & \forall \phi \in \mathbb{V}^{\text{DG}}_1.
  \end{align}
  \end{subequations}
  \endgroup
Finally, we consider a \textbf{DG-primal system} for $T$ in $\mathbb{V}^{\text{DG}}_2$. For this purpose, we use a discretization similar to the conductivity tensor approach in \cite{green2022efficient, vogl2022mesh}, except that we use a split isotropic-anisotropic formulation of the continuous equations in order to have a scheme that is equal to the mixed DG formulation up to the treatment of the purely anisotropic term. The temperature equation is then given by
  \begin{align}
  \left\langle \phi, \pp{T_h}{t} \right\rangle - \text{IP}_b(T_h; \phi) - \text{IP}(T_h; \phi)  = \langle \phi, S\rangle && \forall \phi \in \mathbb{V}^{\text{DG}}_2, \label{primal_DG}
  \end{align}
where the isotropic term $\text{IP}$ is discretized according to \eqref{IP_term}, and the anisotropic term $\text{IP}_b$ is given by \eqref{IP_term_b}. The interior penalty parameter for the resulting isotropic DG-Laplacian is set as for the mixed DG setup, while the one for the anisotropic DG-Laplacian is set to $10$.
As before, the above space discretizations are discretized in time using the implicit midpoint rule. The exception for the inflow boundary term involving $\zeta_\text{in}$ -- as discussed at the end of \Cref{sec_spat_disc} -- does not apply here, since the above space discretizations do not contain such a term.

Next to alternative discretizations, we also consider an alternative solver strategy for the DG and CG mixed systems, which is given by the standard block matrix form \eqref{mat-discr} with the natural (1,1) Schur complement setup as a precondioner to an outer flexible GMRES Krylov iteration. The lower right block, \ccrf{which corresponds to the block diagonal DG auxiliary variable mass matrix $M_\zeta$, is inverted directly}{9}. The inverse Schur complement 
\begin{equation}
\tilde{S}^{-1} = \big(\tfrac{1}{\Delta t}M_T + \kappa_\perp L + \kappa_\Delta G_b^T M_\zeta^{-1}G_b\big)^{-1}, \label{S_standard}
\end{equation}
is solved for using a conjugate gradient Krylov iteration, using BoomerAMG from hypre \cite{falgout2002hypre} as a preconditioner. The multigrid options are left to the default ones. The outer and inner tolerances are set to $10^{-8}$ and $10^{-3}$, respectively.

\subsection{Implementation}\label{Ch3-Implementation} The discretization and solver were implemented using Firedrake ~\cite{Rathgeber2016}, which heavily relies on PETSc \cite{balay2019petsc}. A parallel version of AIR is implemented in hypre \ccrf{(see e.g.~\cite{hanophy2020parallel} for scalability results), but the PETSc interface is incomplete, and we are currently unable to access all necessary features through PETSc and Firedrake.}{ii}. As mentioned in \Cref{sec:solver:solver}, we therefore instead call AIR through PyAMG. The corresponding Python functions are then used in the Firedrake framework via a custom preconditioner interface. The resulting setup can only be used in serial, and the solver efficiency results below therefore show wall-clock times for serial runs.\footnote{\ccrf{We are developing a general and parallel version of block sparse AIR with additional features not available in \emph{hypre} in the AMG library RAPtor \cite{BiOl2017}, which we will use in future work.}{ii}}

The diffusion-based multigrid based solver strategy described in the previous subsection uses a classical AMG solver as implemented in hypre's BoomerAMG \cite{falgout2002hypre}. The {hypre} implementation allows for more options and is more geared towards efficiency than the PyAMG version of classical multigrid; indeed we found the BoomerAMG implementation to generally require slightly fewer iterations, as well as less wall-clock time. Note, we also tested several AMG methods implemented in PyAMG specifically designed for anisotropic diffusion, including \cite{manteuffel2017root,notay2010aggregation,Schroder.2012}, but again found the hypre classical AMG implementation to consistently win in terms of iteration count and wall-clock time.

For the convergence results, all systems were solved with a direct solver. DG based discretizations are solved using \ccrf{MUMPS}{10}\ccbf{}{3} \cite{amestoy2000multifrontal}, as this is less taxing on the memory (DG spaces have more degrees of freedom than the corresponding CG spaces). For CG, \ccrf{MUMPS}{10}\ccbf{}{3} did not consistently converge at higher anisotropy ratios and higher resolutions, and so \ccrf{SuperLU\_DIST}{10} \cite{li2003superlu_dist} was used instead. For problems where both of these choices worked, the results were the same up to round-off; the choice of solver is for computational speedup and does not affect the solution accuracy.

Next to the DG scheme described in \Cref{sec_spat_disc}, as well as the above schemes used for comparison purposes, we also attempted an AIR based strategy for a mixed CG discretization. However, we found such a scheme difficult to implement in view of the strong Dirichlet boundary conditions. In typical finite element codes, it is non-trivial to extract the resulting symmetric transport operators, noting that we need the latter operators for our block-row swapped solver strategy. Other than the implementation difficulty, we further thought it justified to skip the mixed CG discretization together with the AIR framework, since the mixed DG formulation leads to a lower connectivity than the CG one, which is advantageous for the latter framework.

\subsection{Convergence Studies} \label{sec_conv_studies}
Here, we test the convergence order of our novel space discretization, using a variation of a test case considered in \cite{gunter2007finite} and originally presented in the paper introducing the magnetohydrodynamic model NIMROD \cite{sovinec2004nonlinear}. The original test case from the literature is posed in 2 dimensions; for the purposes of our problem, we extrude this into 3D in order to avoid having to deal with the magnetic field singularity (with $|\mathbf{B} | = 0$) present in the 2D problem, which we found to be problematic for the upwinding scheme. Since such singularities do not occur in physically realistic scenarios, we do not consider this as a weakness of our scheme. The domain is described by $\Omega = [0, L_x]\times[0, L_x]\times [0, L_z]$, which is periodic in the $z$-direction. We take $L_x = 1$, $L_z = 5$. The initial temperature field is given by
\begin{equation}
T_0 = \sin\left(\frac{\pi x}{L_x}\right)\sin\left(\frac{\pi y}{L_x}\right), \label{T_ic_NIMROD}
\end{equation}
with $T_\text{BC} = 0$ on $\partial \Omega$. The magnetic field $\mathbf{B}$ is chosen to align with the contours of $T_0$ in the two base dimensions and be constant in the third; namely that
\ccrf{
\begin{equation}
\mathbf{B} = (B_x, B_y, B_z) = (-\partial_y T_0, \partial_x T_0, 5), \label{B_ic_NIMROD}
\end{equation}
}{11}
which ensures that $(B_x, B_y) = \nabla^{\perp}T_0$\, for 2D curl $\nabla^\perp = (-\partial_y, \partial_x)$. \ccbf{Note that the magnetic field line configuration is such that as in the original 2D version, all magnetic field lines are \textit{closed}, in the sense that they do not intersect with the boundary.}{i.}. Finally, as in the original test case, we add a counter forcing $S$ given by $-\kappa_\perp \Delta T_0$, which is included in order to ensure a steady state test case \ccrf{with analytic solution $T(t) = T_0$}{12}. The error can then be measured relative to the initial condition expression \eqref{T_ic_NIMROD} as
\begin{equation}
e_T(t) = \frac{\|T_h(t) - T_0\|_2}{\|T_0\|_2}. \label{L2_error_T}
\end{equation}

As described in \Cref{sec_spat_disc}, in meshing the domain, we use a triangular mesh for $[0, L_x]\times[0, L_x]$, which is then extruded in the third dimension. Since the initial conditions are highly symmetric, we use a perturbed version of a regular triangular mesh. This ensures that our convergence study really is analyzing the effect of the discretization on the error, as opposed to symmetries in the particular chosen setup. \ccrf{The base mesh and initial conditions on the extruded mesh}{14}\ccbf{ are depicted in Figure \ref{fig_base_mesh_conv_ic}}{i.} in \Cref{appendix:base_mesh}. \ccrf{Next to considering possible symmetries, we recall that spurious cross-field pollution largely results from a non-alignment of magnetic field lines with the mesh. The perturbed mesh together with a non-zero $z$-component for $\mathbf{B}$ -- leading to helical rather than circular field lines -- ensures that the 1D magnetic field lines intersect with the 2D mesh facets in a variety of angles. This renders the test case 3D for the purpose of anisotropic diffusion, even if both $T_0$ and $\mathbf{B}$ do not vary in the $z$-component.}{iii.}

Three resolutions are tested, starting from resolution $\Delta x = L_x/7$ for the base mesh, which is then twice refined by splitting cells. For the resolution in the extruded direction, we note that due to the setup of the initial conditions, we find that 2 cells are sufficient to resolve the latter direction, and increasing the number of cells will not decrease the overall error. To keep the overall computational cost low, we therefore use 2 cells in the extruded direction for all three base resolutions. Finally, we use a time step of $\Delta t = 10^{-3}$ for $100$ time steps, giving a final time of $t_{max} = 0.1$. We test anisotropy ratios of $10^3$, $10^6$ and $10^9$ -- implementationally, we set $\kappa_\perp = 1$, such that the anisotropy ratio is equivalent to $\kappa_\parallel$. The resulting convergence plots are presented Figure \ref{order_plot}, where we consider $e_T(t)$ averaged over the last two time steps. This is done since for the lower anisotropy ratio $10^3$, we found all four discretizations to exhibit a small degree of oscillatory behavior in the error across time steps, likely due to an interplay between the isotropic diffusion and counter-forcing terms.
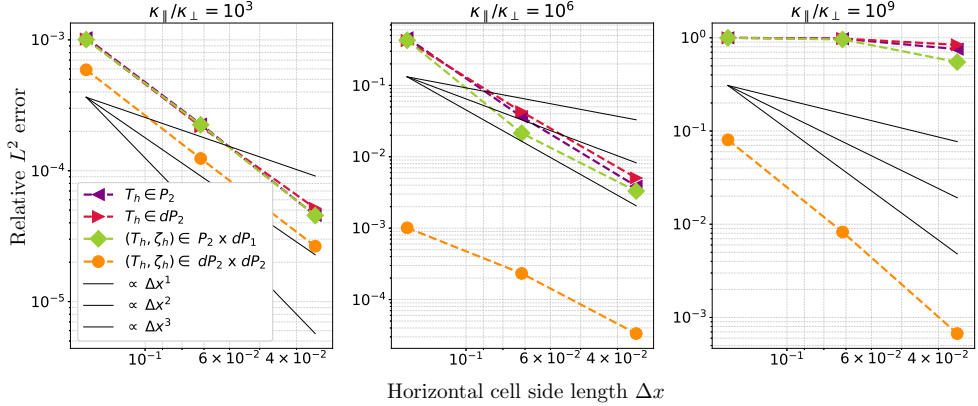
\begin{figure}[ht]
\vspace{-2mm}
\begin{tikzpicture}
  \node[rotate=90,scale=0.8] {\hspace{2cm}Relative $L^2$ error};
\end{tikzpicture}
\hspace{-3mm}
\begin{subfigure}{.32\textwidth}
  \centering
  \includegraphics[width=1.0\linewidth]{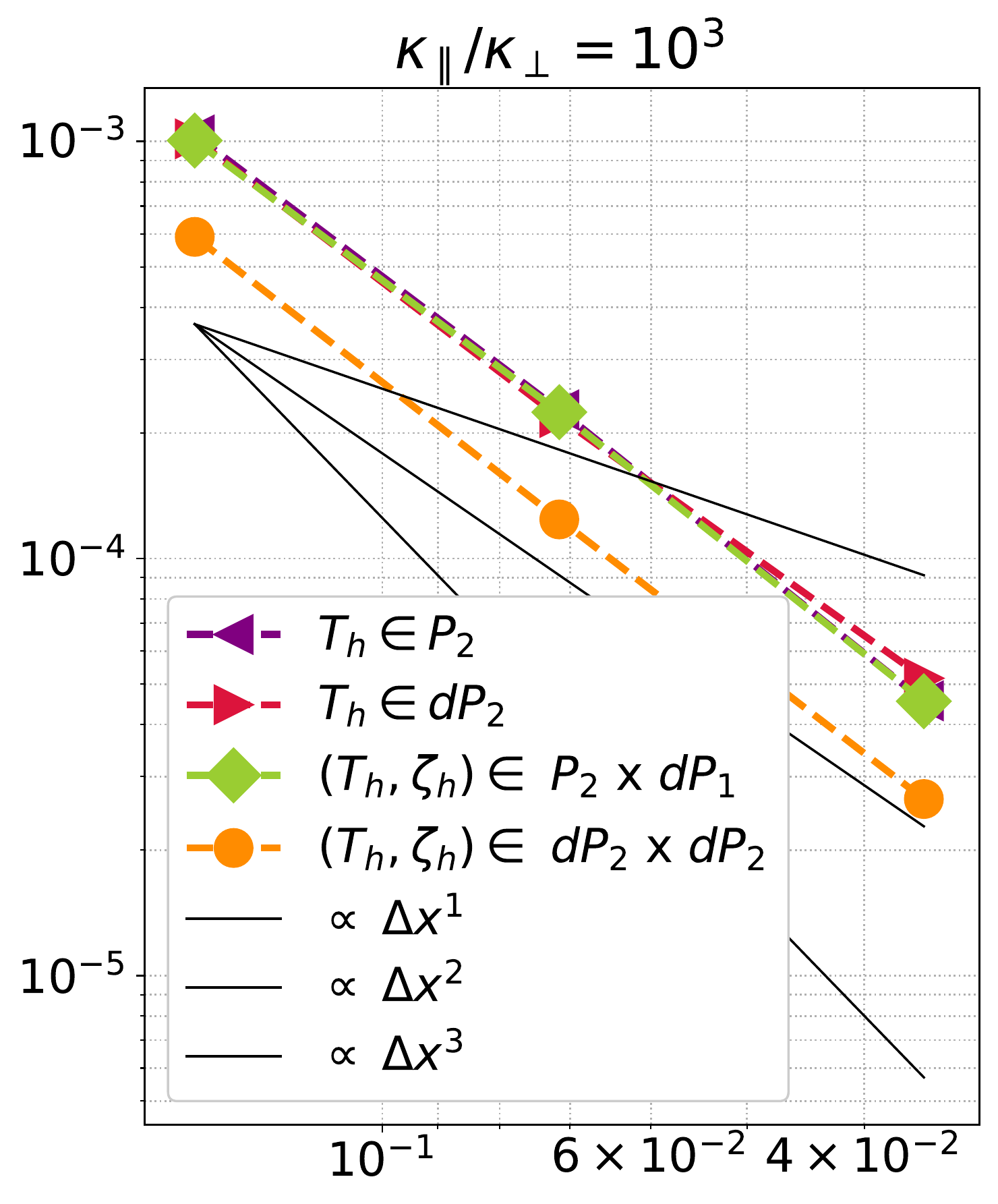}
\end{subfigure}
\begin{subfigure}{.32\textwidth}
  \centering
  \includegraphics[width=1.0\linewidth]{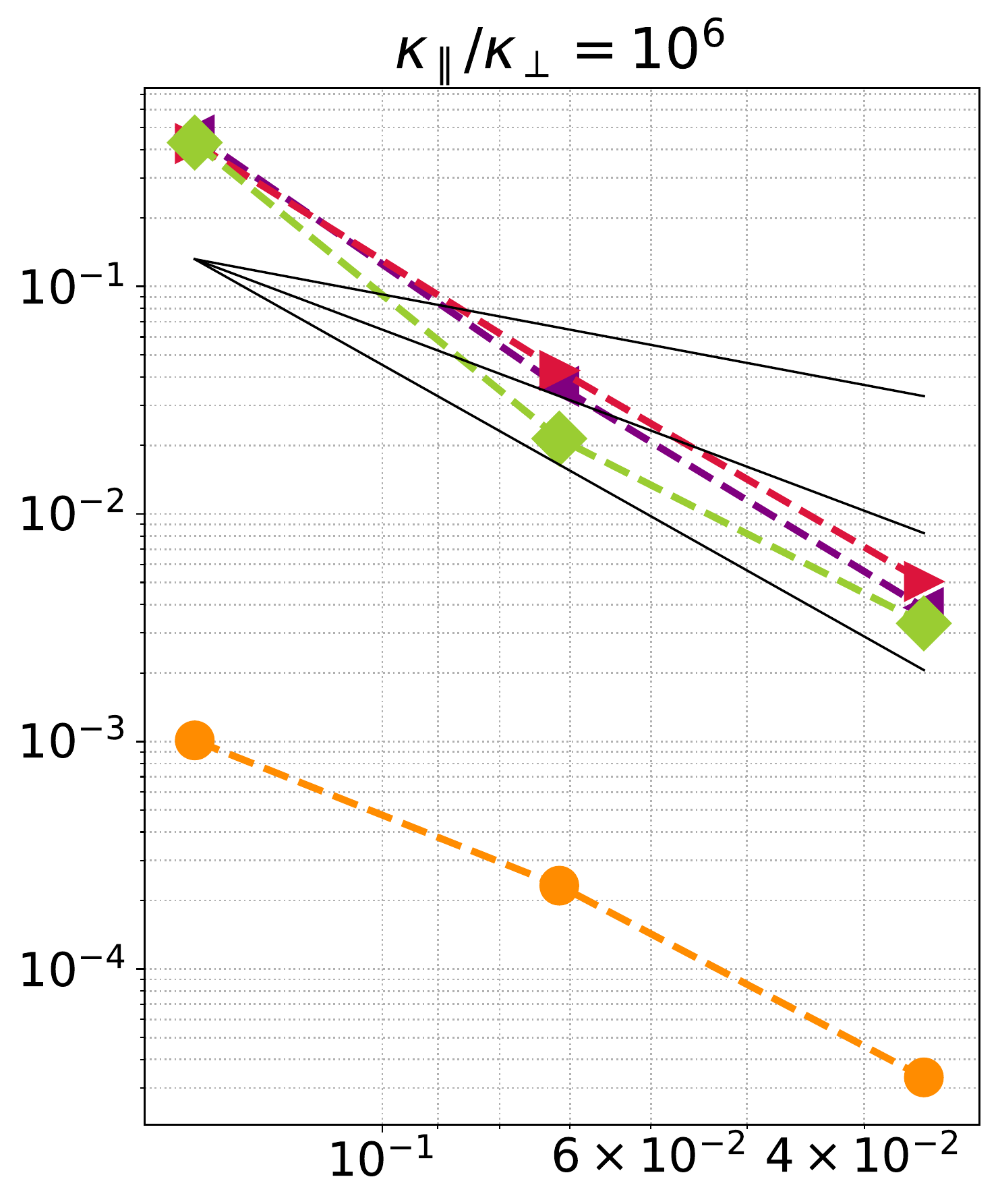}
\end{subfigure}
\begin{subfigure}{.32\textwidth}
  \centering
  \includegraphics[width=1.0\linewidth]{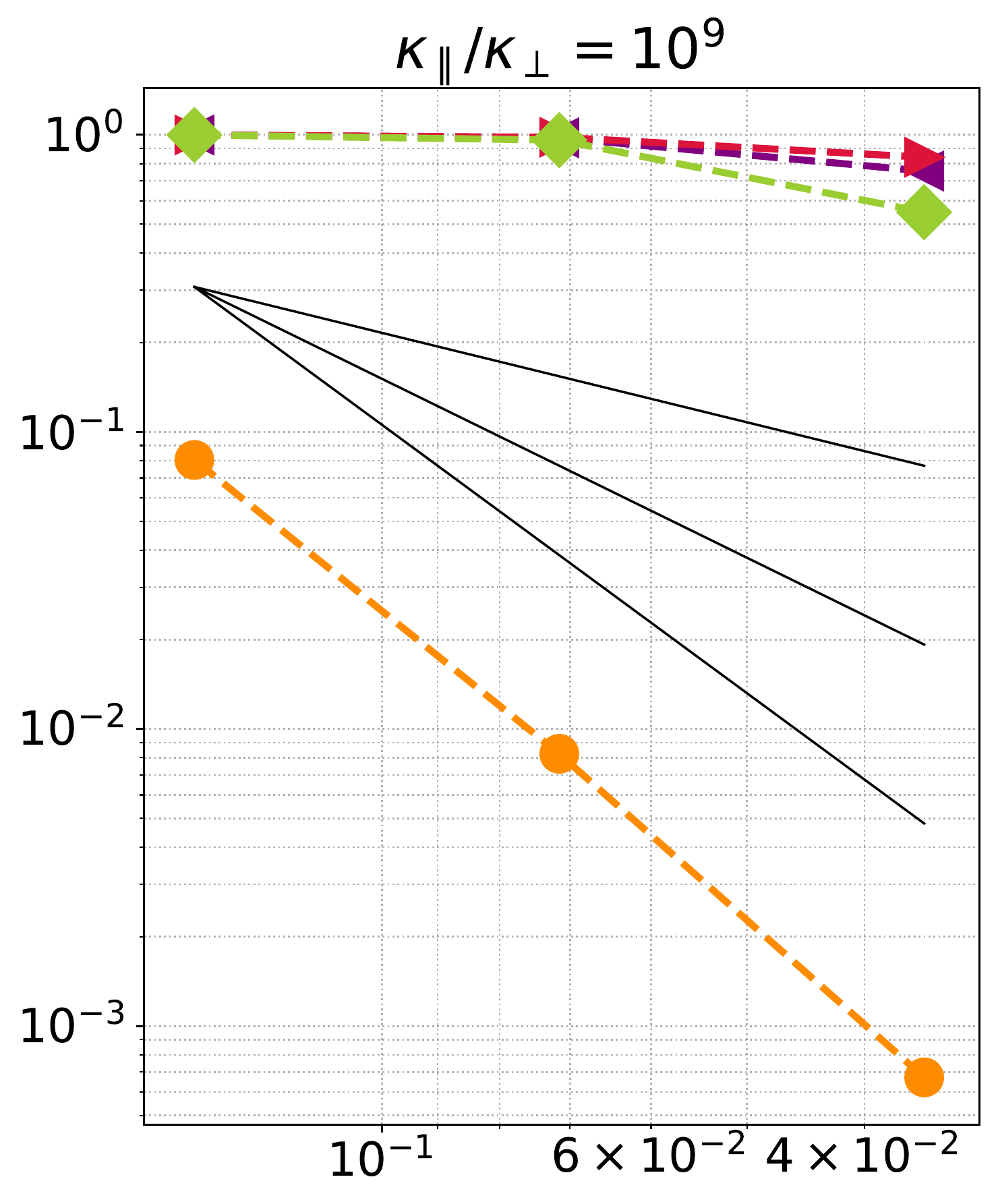}
\end{subfigure}
\vspace{-1mm}
\begin{tikzpicture}
  \node[scale=0.8] {\hspace{6.3cm}Horizontal cell side length $\Delta x$};
\end{tikzpicture}
\vspace{-1mm}
\caption{Relative $L^2$ error \eqref{L2_error_T} convergence plots for the two non-mixed and mixed discretizations, respectively. The anisotropy ratios are given from left to right by $10^3$, $10^6$, and $10^9$. Orange circles denote the novel mixed DG scheme \eqref{DG_upw_discr}, green diamonds the mixed CG scheme \eqref{mixed_CG}, red triangles the primal DG scheme \eqref{primal_DG}, and purple triangles the primal CG scheme \eqref{primal_CG}. The black lines indicate first to third order convergence.} \label{order_plot}
\end{figure}

For an anisotropy ratio of $10^3$, all four schemes show similar errors and convergence rates. For $10^6$, all schemes still display second to third order convergence rates; however, the errors for the newly introduced scheme \eqref{DG_upw_discr} are approximately two orders of magnitude smaller than for the other three schemes. Finally, for the high degree of anisotropy $10^9$, scheme \eqref{DG_upw_discr} still displays favorable error values and convergence rates, while the remaining schemes exhibit relative errors of order $1$. Inspecting the runs' field development, we found these large errors likely to be due to a large spurious contribution of the heat flux' parallel component in the perpendicular direction, thereby rendering the latter direction overly diffusive. The mixed DG scheme, on the other hand, is able to avoid such an excessive spurious diffusion likely due to a better representation of cross-cell fluxes through the upwind formulation. Note that the results for the mixed CG scheme are considerably worse than in \cite{gunter2007finite}. We found this to be unrelated to using a 3D version of the original test case; instead, the difference seems to be related to the use of triangles in the base mesh in this work, versus quadrilaterals in \cite{gunter2007finite}. Switching to  \ccrf{a quadrilateral base mesh perturbed analogously to the triangular one}{15}, in additional tests not shown here, we found the performance of the mixed CG scheme to be drastically improved; albeit with an error that was still larger than the one for the mixed DG scheme.

Finally, we again stress a caveat of the newly introduced scheme \eqref{DG_upw_discr}, given by the heat flux guess $\zeta_\text{in}$ at the inflow boundary. While non-dynamical test cases -- such as the one considered here -- allow for an analytic solution given by the initial condition, they may conceal problems related to dynamics that vary in time. That said, in further tests with quadrilateral meshes, a non-steady dynamical evolution, and an anisotropy ratio of $10^3$, we found little difference between the mixed CG and the mixed DG discretizations.

\subsection{Solver Studies}\label{sec_efficiency}
While we used a modification of the test case from \cite{sovinec2004nonlinear} for the convergence results, here we consider a new test case, which in particular has open field lines only. We consider the same base mesh and domain as in the previous test case. Again, the base mesh starts with resolution $\Delta x = L_x/7$, and we use $2$ cells in the extruded direction. As before, we refine twice by splitting cells. However, this time, we also refine the number of cells in the extruded direction in order to ensure a truly three-dimensional solver study. For the highest refinement level, this leads to approximately $10^5$ and $4.5\times10^5$ degrees of freedom for $\mathbb{V}^\text{CG}_2$ and $\mathbb{V}^\text{DG}_2$, respectively. The initial temperature field is given by
\begin{equation}
T_0 = 1 + \frac{1}{20}\left(1 - \cos\left(\frac{2 \pi y}{L_x}\right) \right)\sin\left(\frac{\pi x}{L_x}\right) + x + \frac{y}{10}, \label{T_ic_conv}
\end{equation}
with $T_\text{BC}$ defined according to the values of $T_0$ on $\partial \Omega$. As in the convergence study, the magnetic field $\mathbf{B}$ is chosen to align with the contours of $T_0$ in the two base dimensions. This leads to field lines which are slightly slanted with respect to the domain, as well as curved towards the domain's center. Further, $\mathbf{B}$ is again chosen to be constant in the third dimension; this time, we set $B_z = \tfrac{15}{2}$. \ccrf{As for the convergence study, we note that the test case is 3D for the purpose of anisotropic diffusion, as this setup ensures that the 1D magnetic field lines intersect with the 2D mesh facets in a variety of angles. The effect of this versus a 2D mesh can for instance be seen in the largest eigenvalue analysis as described at the end of \Cref{sec:solvers:eigs}.}{iii.} Finally, since here we only consider the solver performance and not convergence rates, there is no need for an analytic solution and the associated counter-forcing; hence we set $S=0$.

\ccbf{We begin by demonstrating the need for multilevel solvers versus simple approximate inverse preconditioners like ILU, \emph{independent of the parallel benefits of multigrid over coarse approximate inverses.} \Cref{ILU-table} shows iterations to convergence for GMRES preconditioned with ILU(2), applied to the anisotropic inner Schur complement solve \eqref{S_standard} for the mixed CG scheme \eqref{mixed_CG}, for three anisotropy ratios, with $\Delta t=10^{-3}$, and the medium spatial refinement level used in this section. For moderate anisotropies, ILU(2) performs well, but at $\kappa_\|/\kappa_\perp = 10^7$, after 30,000 ILU(2) iterations the residual has decreased by less than a factor of $10^{-3}$. Further accounting for the poor scaling of ILU and many basic sparse approximate inverses in problem size and particularly number of processors, we see that for realistic problems and/or large-scale simulation, more advanced preconditioners are \emph{necessary.}
\begin{table}
\begin{center}
{\setlength{\extrarowheight}{3pt}
\setlength{\tabcolsep}{3.7pt}
\begin{tabular}{|c||c|c||c|c||c|c|}
 \hline
$\kappa_\|/\kappa_\perp$ & \multicolumn{2}{c||}{$10^4$} & \multicolumn{2}{c||}{$10^5$} & \multicolumn{2}{c|}{$10^7$}\\\hline
Outer its. & Inner its. & Outer res. & Inner its. & Outer res. & Inner its. & Outer res. \\
 \hline
1&20&2.7e-5&119&5.1e-4&10000&9.97e-1\\
2&21&1.8e-11&10000&4.2e-7&10000&1.5e-2\\
3&&&10000&2e-9&10000&5.8e-3\\
  \hline
\end{tabular}}
\caption{Convergence results for GMRES preconditioned by ILU(2) applied to the mixed CG scheme \eqref{mixed_CG} for three anisotropy ratios. Inner iterations are solved to $10^{-6}$ relative residual, and outer stopping tolerance is $10^{-8}$ relative residual.} \label{ILU-table}
\end{center}
\end{table}
}{iii.}

Now we move on to consider multigrid solvers, and our novel advection-based preconditioner. We consider results for a) the mixed DG-scheme \eqref{DG_upw_discr} together with the AIR solver strategy described in \Cref{sec:solver:solver}, b) the mixed DG-scheme \eqref{DG_upw_discr} together with the classical AMG based solver strategy described in \Cref{comp_disc_solver}, as well as c) the mixed CG scheme \eqref{mixed_CG} together with the aforementioned classical AMG based solver strategy. As before, the time step is set to $\Delta t = 10^{-3}$, and we consider iteration counts as well as wall-clock time from the $2^{nd}$ to the $5^{th}$ time step, \ccrf{i.e.~}{1}a total of four time steps. Note that we disregard the first time step as a simple way of avoiding to count wall-clock time associated with form assembly. Further, runs which either exceed a specified average wall-clock time or outer iteration count per time step (approx. $1500$ seconds and $10^4$ iterations, respectively) are discarded. The resulting averages per time step for the total inner iteration counts, \ccbf{outer iteration counts}{iii.}\ccrf{}{8, 17}, as well as wall-clock times, for anisotropy ratio $\kappa_\parallel/\kappa_\perp = 10^k$, $k = 2, ..., 10$, are depicted in Figures \ref{Efficiency_iterations}, \ref{Efficiency_outer_iterations}, and \ref{Efficiency_wallclock}, respectively.

\begin{figure}[ht]
\begin{center}
\vspace{-2mm}
\begin{tikzpicture}
  \node[rotate=90,scale=0.8] {\hspace{3cm}Avg. total inner iteration count per time step};
\end{tikzpicture}
\hspace{-3mm}
\includegraphics[width=0.95\textwidth]{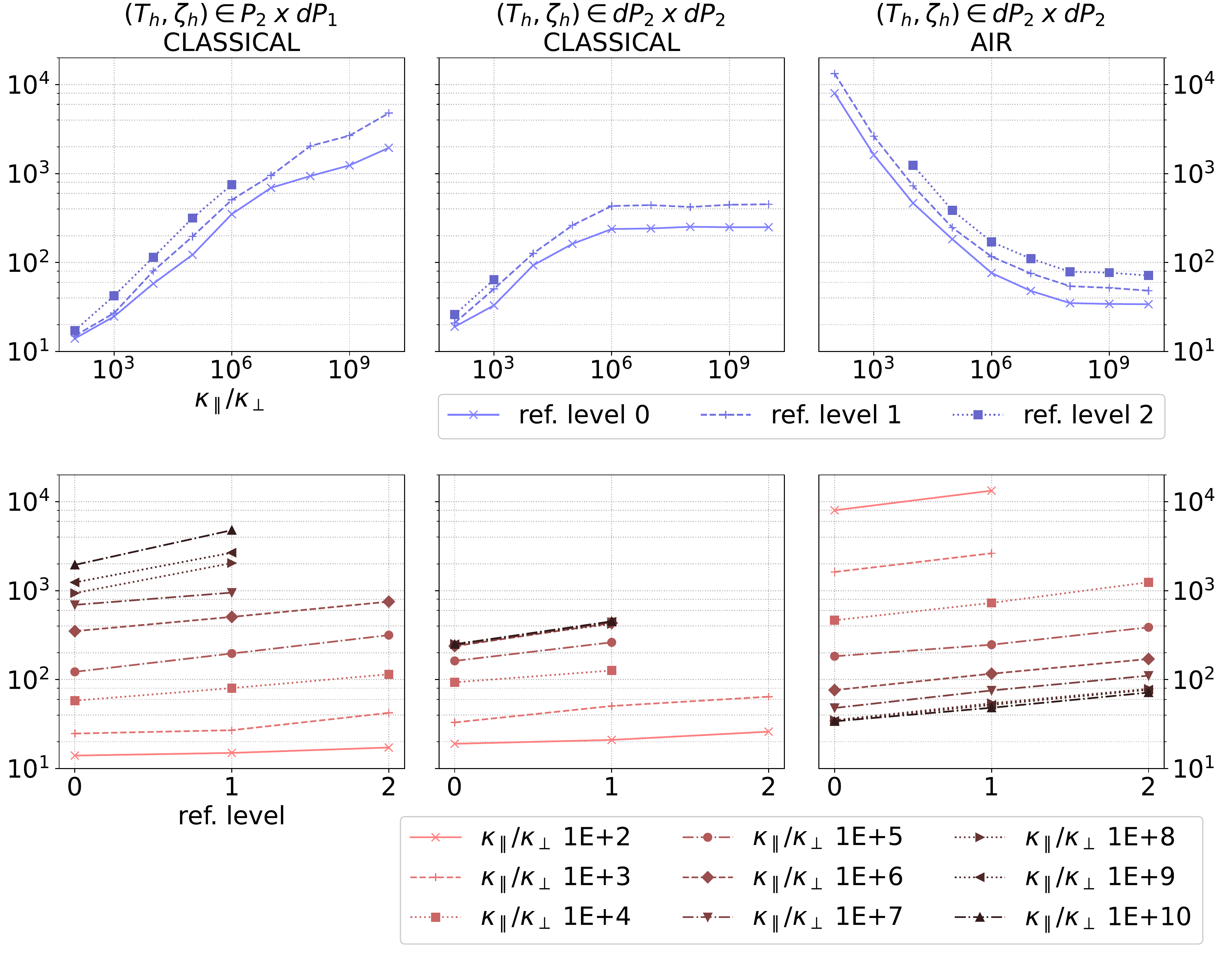}
\vspace{-1mm}
\caption{Average total inner iteration counts per time step using setup described in \Cref{sec_efficiency}. Top row: counts with respect to anisotropy ratio. Bottom row: counts with respect to refinement level. Left column: mixed CG scheme \ccrf{with inner iterations corresponding to solving for \eqref{S_standard} using classical AMG. Center column: mixed DG scheme with inner iterations corresponding to solving for \eqref{S_standard} using classical AMG. Right column:  mixed DG scheme with each inner iteration corresponding to solving for either \eqref{inner_solve_AIR_inv} or \eqref{inner_solve_AIR_inv_T}, using AIR}{16}.}\label{Efficiency_iterations}
\end{center}
\end{figure}

\begin{figure}[ht]
\begin{center}
\vspace{-2mm}
\begin{tikzpicture}
  \node[rotate=90,scale=0.8] {\hspace{3cm}Avg. outer iteration count per time step};
\end{tikzpicture}
\hspace{-3mm}
\includegraphics[width=0.95\textwidth]{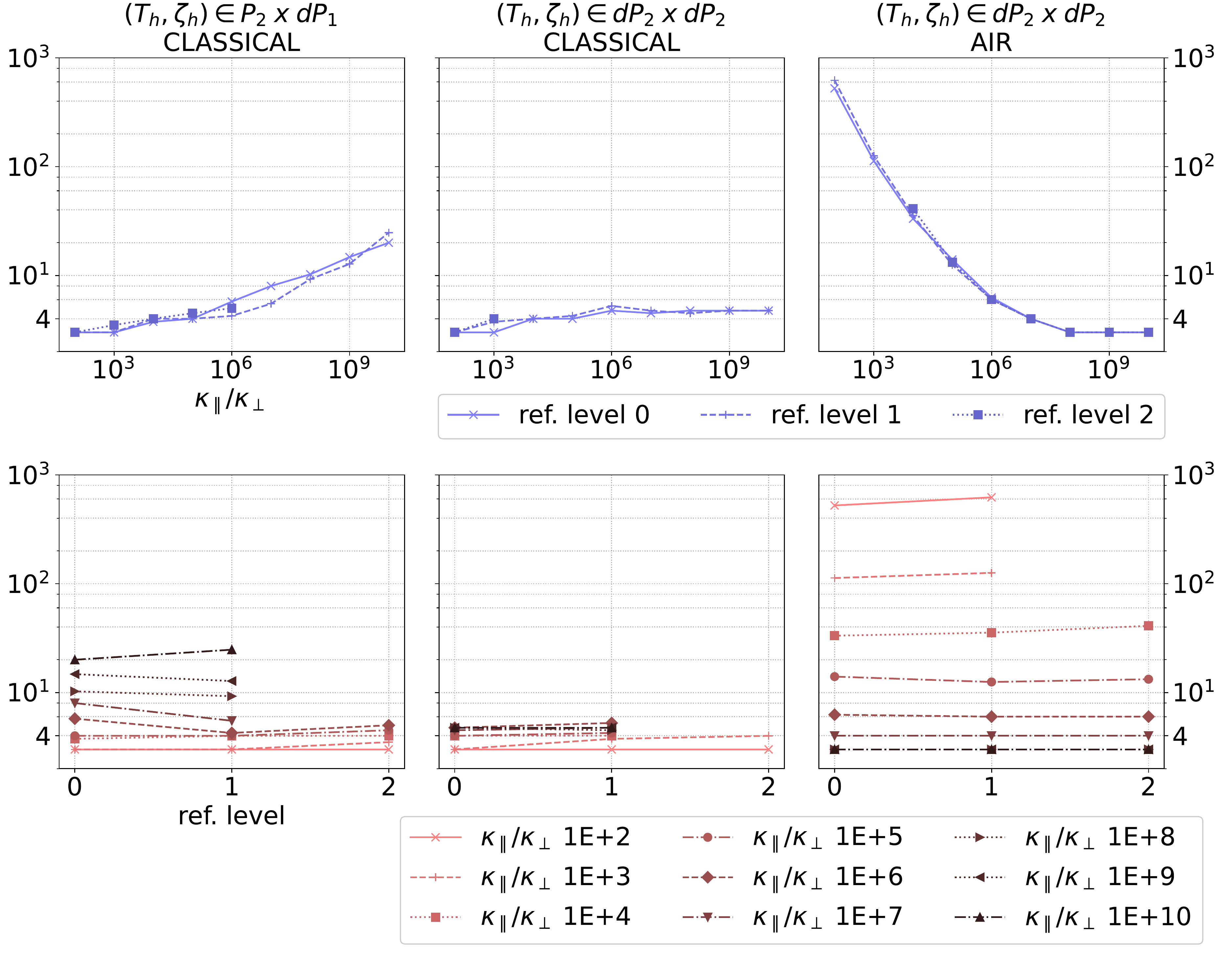}
\vspace{-1mm}
\caption{\ccbf{Average outer iteration counts per time step using setup described in \Cref{sec_efficiency}, with inner residual tolerance $10^{-3}$. Top row: counts with respect to anisotropy ratio. Bottom row: counts with respect to refinement level. Left column: mixed CG scheme with inner iterations corresponding to solving for \eqref{S_standard} using classical AMG. Center column: mixed DG scheme with inner iterations corresponding to solving for \eqref{S_standard} using classical AMG. Right column:  mixed DG scheme with each inner iteration corresponding to solving for either \eqref{inner_solve_AIR_inv} or \eqref{inner_solve_AIR_inv_T}, using AIR.}{iii.}\ccrf{}{8, 17}}\label{Efficiency_outer_iterations}
\end{center}
\end{figure}

\begin{figure}[ht]
\begin{center}
\begin{tikzpicture}
  \node[rotate=90,scale=0.8] {\hspace{4cm}Avg. wall-clock times per time step};
\end{tikzpicture}
\includegraphics[width=0.95\textwidth]{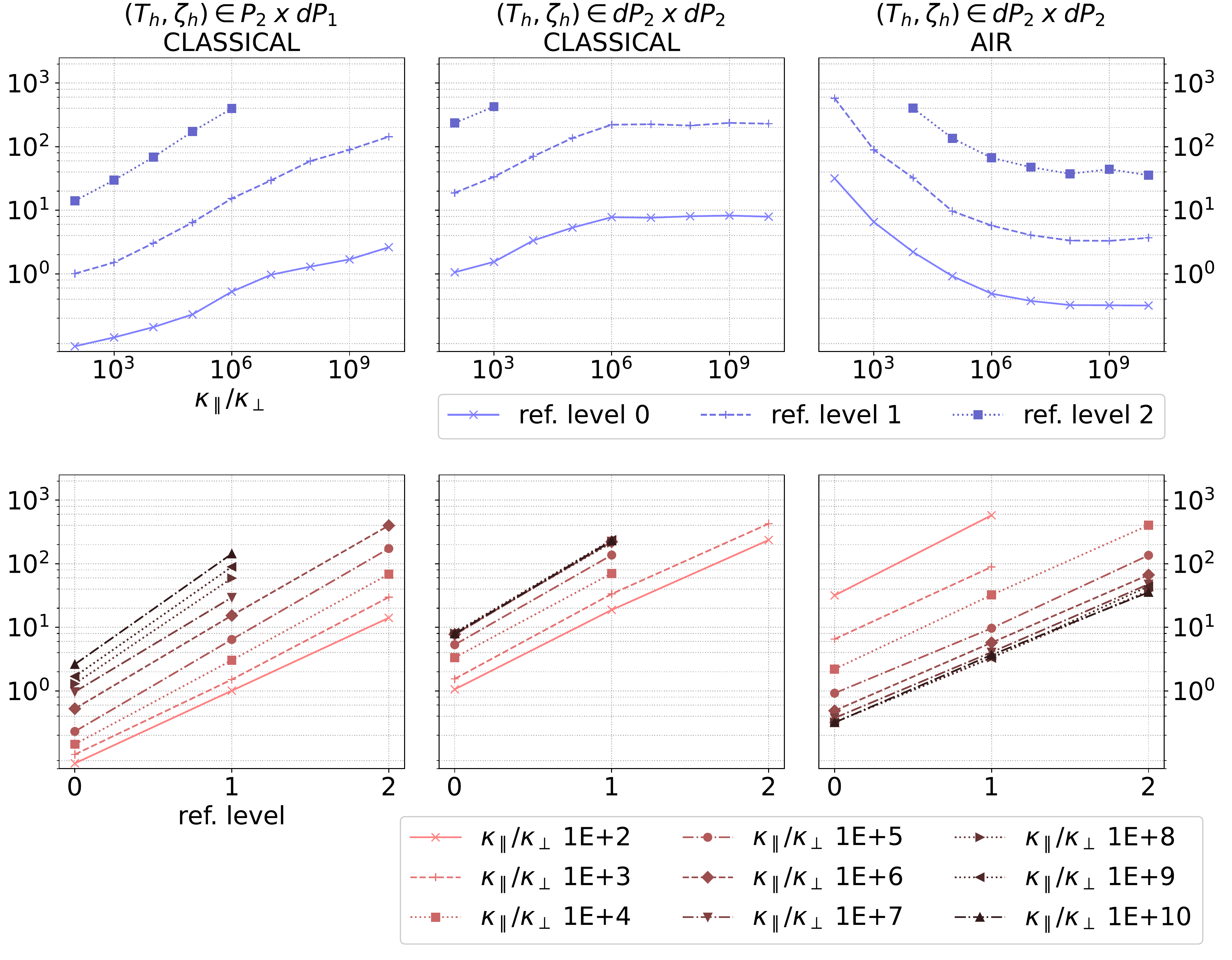}
\vspace{-1mm}
\caption{Average wall-clock times in seconds per time step using setup described in \Cref{sec_efficiency}. Top row: times with respect to anisotropy ratio. Bottom row: times with respect to refinement level. Left column: mixed CG scheme\ccrf{, solving for \eqref{S_standard} using classical AMG. Center column: mixed DG scheme, solving for \eqref{S_standard} using classical AMG. Right column:  mixed DG scheme, solving for \eqref{inner_solve_AIR}, using AIR.}{16}}\label{Efficiency_wallclock}
\end{center}
\end{figure}

As expected for the mixed CG scheme solved using classical AMG, we find a steady increase in solver iteration counts and according wall-clock times as the ratio of anisotropy increases. 
For the highest refinement level, runs with anisotropy ratio $10^7$ or higher did not meet the imposed cutoff criteria. Passing from the first to the second refinement level, runs at the highest three anisotropy ratios experience a wall-clock time increase by a factor of approximately $100$. Surprisingly, for the mixed DG scheme solved using classical AMG, we find the increase in inner solver iteration count and wall-clock times to eventually stall as we increase the ratio of anisotropy. For the lowest refinement level, the times stall at approximately $8$ seconds, while for the next higher one, they do at about $220$ seconds -- which indicates a roughly 27-fold increase across the two refinement levels. For this setup, most of the highest refinement level runs did not meet the imposed time limit. In comparison to the mixed CG scheme solved using classical AMG, this is likely due to the much higher number of degrees of freedom associated with the DG scheme.

For the mixed DG scheme solved using our newly proposed solver strategy including AIR, we instead find a \textit{decrease} in inner iteration counts and wall-clock times as the anisotropy ratio increases. As described in \Cref{sec:solvers}, this is as expected, as such higher ratios fit better into our inner solver strategy's favorable parameter regime \eqref{parameter_regimes}. 
For the highest refinement level, runs with anisotropy ratio $10^3$ or lower did not meet the imposed time limit. In particular, the advection-based block preconditioner struggles with lower anisotropy ratios, while it performs very well for ratios $10^8$ and higher. For instance, at the highest anisotropy ratio $10^{10}$ and the medium refinement level, the mixed CG formulation solved with classical AMG requires an average of approximately $5000$ iterations, while the mixed DG formulation solved with our new AIR based strategy only requires approximately $50$, \ccrf{i.e.~}{1}around $100$ times fewer. Given the aforementioned rates of increase for the number of iterations as we increase the refinement level, we expect this gap to rapidly increase for higher resolutions. Similarly, for the same anisotropy ratio and refinement level, we find the mixed DG/AIR based run to be around $35$ times faster than the mixed CG/classical AMG based run; this holds true even though the DG discretization has far more degrees of freedom than the CG one, and we use a classical AMG based code implementation that is more streamlined for efficiency than our AIR based one. Again, given the rates of increase in wall-clock time as we increase the refinement level, we expect this gap to be larger still for higher resolutions, provided that the classical AMG based approach still converges. In addition, recent work has proposed a modified AIR \cite{new-air} that reduces total computational time compared with $\ell$AIR as implemented in hypre \cite{manteuffel2019nonsymmetric,manteuffel2018nonsymmetric} by several times, suggesting further multiplicative speedups possible.

Finally, we find the difference in iteration count for the mixed DG scheme at the three highest anisotropy ratios to not be as pronounced between the AIR- and classical AMG based solver strategies. Picking again the example of an anisotropy ratio of $10^{10}$ at the medium refinement level, we find that the AIR based strategy requires about $10$ times fewer iterations than the classical AMG based one. It is unclear why the classical AMG based strategy performs this well in terms of iteration count for the mixed DG scheme. However, considering the wall-clock time, we find that the mixed DG/classical AMG based run is around $60$ times slower than the AIR based one, with significant degradation as the mesh is further refined.

We end the discussion with observations on the AIR strategy's solver behavior. First, we found the number of inner iterations to remain very confined throughout the range of anisotropy ratios that we tested in this section. For the highest refinement level and an anisotropy ratio of $10^2$, the inner solves for the transport operator and Schur complement take approximately 10 and 22 iterations per outer iteration, respectively, with the larger number of Schur complement iterations due to the modified inner tolerance proposed in \Cref{sec:solver:solver}. For $10^6$, we obtain 10 and 19, respectively. Finally, for $10^{10}$, the counts are 10 and 14, respectively. In other words, the inner solver count for the transport operator remains the same, while for the Schur complement, it slowly decreases. This is because as $\kappa_\perp L$ becomes a more significant off-diagonal contribution to the Schur complement's right-hand side as $\kappa_\Delta$ decreases, we require a larger residual reduction to \ccrf{meet}{19}\ccbf{}{4} an absolute residual tolerance.

Second, considering the relatively low inner iteration counts throughout the range of anisotropy ratios, we note that it is the outer solver that begins to struggle at anisotropy ratios of $10^5$ or less (see \Cref{Efficiency_outer_iterations}). Recall that the Schur complement is preconditioned using the transport operator (see \eqref{eq:prec-S}), which we expect to work well when we are within the valid parameter regime \eqref{parameter_regimes} for our solver. The latter regime increasingly ceases to be satisfied for ratios $10^5$ or less, as can be inferred from the largest eigenvalue \Cref{Eigen_table_3D}. For a ratio $10^5$, the largest eigenvalues of the preconditioned Schur complement are $7.6$, $5$ and $4.4$ for the three spatial refinements, respectively. 
They then further increase inversely proportionally to the anisotropy ratio. Consistent with the spectral analysis (see \Cref{Eigen_table}), the proposed method is efficient in 2D for smaller anisotropies, on the order of $10^3$, but results are not presented due to space, and because we are interested in 3D-simulations in practice. \ccrf{In \Cref{Efficiency_outer_iterations} we see that outer iteration counts are largely scalable in mesh resolution, and when increasing anisotropy ratio from $10^2$ to $10^5$ or higher, the outer iteration count drops from 1000 to $\mathcal{O}(1)$.}{8, 17}\ccbf{}{iii., iv.}

\ccbf{Third, the average number of total inner iterations per time step increases only moderately as the spatial resolution is doubled. For instance, for an anisotropy ratio of $10^8$, the count is approximately $35$, $50$, and $80$, respectively, for the three spatial refinements. This is actually similar to AIR applied to advection, which tends to see slow logarithic growth in iteration count with problem size \cite{manteuffel2019nonsymmetric}, unlike the perfect scaling that can be achieved for isotropic diffusion.}{iii., iv.} Altogether, for the very high anisotropy ratios of $10^6$ or more typically considered in magnetic confinement fusion, our novel mixed DG/AIR based scheme performs very well, and clearly outperforms the classical AMG based approaches.

\section{Conclusion} \label{sec_conclusion}
Confined plasma simulations may exhibit high anisotropy, which needs to be resolved accurately in order to avoid cross-contamination of heat flux from components parallel to the magnetic field into components perpendicular to the magnetic field. This cross-contamination in turn can lead to spuriously short simulated magnetic confinement times. In addition, such heat flux requires implicit integration, and for large-scale simulation requires efficient iterative solvers. In this work, we presented a novel DG space discretization based on the heat flux as an auxiliary variable, and two upwind transport operators. It is coupled to a novel solver strategy based on the AIR AMG method, where the block matrix system resulting from an implicit time discretization is block-row swapped, and is solved for using a Schur complement approach. For high ratios of anisotropy, the corresponding $A_{11}$ block and Schur complement are then both akin to DG upwind transport operators, which can be solved for efficiently using AIR.


We compared the proposed spatial discretization to a CG and a DG direct discretization, as well as another CG mixed discretization which also uses the heat flux as an auxiliary variable. Further, we compared the AIR based solver strategy with a standard one based on classical AMG. In the first test case based on (extruded) prism elements, we showed that for moderate rates of anisotropy, our novel space discretization has comparable errors and an equal convergence rate when compared to the three other discretizations. Further, for higher rates of anisotropy, our discretization retains the low errors and favorable convergence rate (third-order convergence for 2nd-order basis functions), while the remaining discretizations' errors and rates deteriorate. For anisotropy ratios of $10^9$, our proposed discretization achieves error $\approx 1000$ times smaller than the other discretizations tested. In the second set of tests, we demonstrated the newly proposed AIR based solver strategy to be highly efficient at anisotropy ratios $10^6$ and higher. We found the latter strategy to take up to $100$ times fewer solver iterations and $35$ times less wall-clock time than the classical AMG approach at a moderate space resolution, in spite of using a more streamlined and efficient classical AMG implementation than the one for AIR. Additionally, we found the AIR based solver strategy to run within our imposed time and iteration count cut-off criteria at higher resolutions, where this was no longer the case for the AMG approach.

A major shortcoming of the proposed solver methodology -- but not the space discretization -- is that it requires all magnetic field lines to be open. Ongoing work is focused on developing a generalized solver strategy for the proposed discretization that is robust for closed field lines.


%
%
%
%
%

\appendix

\ccbf{\section{Discretization for sheath boundary conditions}\label{appendix:Neumann_BC}
For Neumann type sheath boundary conditions, we separately specify a flux parallel and perpendicular to the magnetic field according to
\begingroup
\begin{subequations} \label{sheath_BC}
\begin{align}
&\mathbf{n}\cdot\mathbf{q}_\parallel = \mathbf{n} \cdot(\kappa_\parallel\nabla_\parallel T) = q_{\parallel,\text{BC}}, \\
&\mathbf{n} \cdot \mathbf{q}_\perp = \mathbf{n}\cdot(\kappa_\perp\nabla_\perp T) = q_{\perp,\text{BC}}, 
\end{align}
\end{subequations}
\endgroup
on $\partial \Omega$. In magnetic confinement fusion applications, the perpendicular heat flux boundary value is typically set to zero, while the parallel heat flux one may be large and may depend on the local temperature, density, and flow speed through the Bohm criterion \cite{schneider2006plasma}.}{ii.}

\ccbf{The Dirichlet boundary condition-based space discretization \eqref{DG_upw_discr} can readily be adjusted to sheath boundary conditions \eqref{sheath_BC}, since the underlying auxiliary variable \eqref{T_eqn_full_zeta} is analogous to the parallel heat flux $\mathbf{q}_\parallel$ occurring in \eqref{sheath_BC}, up to the scalar factor of $\sqrt{\kappa_\Delta}$ versus $(\mathbf{n}\cdot\mathbf{b})\kappa_\parallel$. In the Dirichlet case, we used the known value $T_{\text{BC}}$ and guess $\zeta_{\text{in}}$ in the transport operators \eqref{transport_T_zeta}. For the Neumann type sheath boundary condition case, we swap the roles of known value and guess, yielding
\begingroup
\addtolength{\jot}{2mm}
\begin{subequations} \label{transport_T_zeta_N}
\begin{align}
&L_{b, T,N}(\zeta_h; \phi) \coloneqq L_b(\zeta_h; \phi) + \int_{\partial \Omega_{\text{in}}} \phi \kappa_\Delta \kappa_\parallel^{-1}q_{\parallel,\text{BC}}\;\text{d} S & \forall \phi \in \mathbb{V}^{\text{DG}}_k, \label{transport_T_N}\\
&L_{b, \zeta,N}(\psi; T_h) \coloneqq L_b(\psi; T_h) - \int_{\partial \Omega_{\text{out}}} \sqrt{\kappa_\Delta} T_{\text{in}} (\vect{b} \cdot \mathbf{n}) \psi \;\text{d} S & \forall \psi \in \mathbb{V}^{\text{DG}}_k, \label{transport_zeta_N}
\end{align}
\end{subequations}
\endgroup
where we recall $L_b$ is given by \eqref{transport_out}. In analogy to $\zeta_{\text{in}}$ for the Dirichlet case, here $T_{\text{in}}$ is a value that needs to be estimated; for details, see the paragraph before \eqref{mat-discr-bdry-orig}. To complete the discretization, we include a penalty term for the boundary conditions, analogously to the one for $T_{BC}$ in the Dirichlet case. The space discretization is then given by
\begingroup
\addtolength{\jot}{4mm}
\begin{subequations} \label{DG_upw_discr_N}
\begin{align}
&\left\langle \phi, \pp{T_h}{t} \right \rangle - L_{b, T,N}(\zeta_h; \phi) - \text{IP}_{\!N}(T_h; \phi) = \langle \phi, S \rangle &\text{\qquad $\forall \phi \in \mathbb{V}^{\text{DG}}_k$}, \label{DG_upw_discr_T_N} \\
&\left\langle \psi, \zeta_h \right \rangle \!+\! L_{b, \zeta,N}(\psi; T_h) = -\!\!\! \int_{\partial \Omega}\!\!\!\!\!\! \kappa_{\text{BC}} \psi \big((\mathbf{b}\cdot\mathbf{n})\zeta \!-\! \sqrt{\kappa}_\Delta \kappa_\parallel^{-1}q_{\parallel,\text{BC}}\big) \text{d}S\!\!\!\!\!\!\!\!\!\!\! &\text{\qquad $\forall \psi \in \mathbb{V}^{\text{DG}}_k$},\label{DG_upw_discr_zeta_N}
\end{align}
\end{subequations}
\endgroup
for $\kappa_{\text{BC}} = \tilde{\kappa}_{\text{BC}}h_e$, and where $\tilde{\kappa}_{\text{BC}}$ is non-dimensional. Further, the Neumann type DG Laplacian is given by
\begingroup
\addtolength{\jot}{4mm}
\begin{align}
\begin{split}
\text{IP}_N(T_h; \phi) = &- \left \langle \nabla \phi, \kappa_\perp \nabla T_h \right \rangle
+\int_\Gamma \left \llbracket T_h \right \rrbracket \left \lbrace \kappa_\perp \nabla \phi \right \rbrace \;\text{d}S
+ \int_\Gamma \left \llbracket \phi \right \rrbracket \left \lbrace \kappa_\perp \nabla T_h \right \rbrace \;\text{d}S\\
&\hspace{-12mm} +\!\! \int_{\partial \Omega}\!\! (\kappa_\perp\kappa_\parallel^{-1}q_{\parallel,\text{BC}} + q_{\perp,\text{BC}}) \phi \;\text{d}S \!-\!\! \int_\Gamma \frac{\kappa_\perp \kappa_p}{h_e} \left \llbracket \phi \right \rrbracket \left \llbracket T_h \right \rrbracket \;\text{d}S, \label{IP_term_N}
\end{split}
\end{align}
\endgroup
noting that the sum of parallel and perpendicular heat fluxes (second to last term in \eqref{IP_term_N}) corresponds to the total flux $\kappa_\perp \mathbf{n} \cdot \nabla T$.}{ii.}

\ccbf{The discretization for sheath boundary conditions can be shown to lead to structural statements analogous to those in \Cref{Remark:structure}, with $\sqrt{\kappa}(\mathbf{b}\cdot \mathbf{n})\zeta_{\text{in}}$, $T_{\text{BC}}$ replaced by $q_{\parallel,\text{BC}}$, $T_{\text{in}}$, respectively, and a boundary flux based penalty instead of a boundary temperature value based one in \eqref{discrete_diffusion} and \eqref{energy_discr}. Finally, after combining with an implicit midpoint rule in time, the resulting algebraic system reads
\begin{equation}\label{mat-discr-bdry-N}
\begin{pmatrix}
  \frac{1}{\Delta t}M + \kappa_\perp L_N && \sqrt{\kappa_\Delta} G_b^T \\
  -\sqrt{\kappa_\Delta}G_b && M + \tilde{\kappa}_\text{BC} M_{\text{BC}, h_e}
\end{pmatrix}
\begin{pmatrix}
  T_h^{n+1} \\
  \zeta_h^{n+1}
\end{pmatrix}
=
\begin{pmatrix}
  F_{T,N} \\
  F_{\zeta,N}
\end{pmatrix},
\end{equation}
with notation as for the Dirichlet boundary condition-based system \eqref{mat-discr-bdry}, and where $L_N$, $F_{T,N}$, $F_{\zeta,N}$ denote the DG Laplacian as well as right-hand sides each adjusted for Neumann boundary terms. In particular, the Dirichlet and Neumann type $2\times 2$ block systems \eqref{mat-discr-bdry} and \eqref{mat-discr-bdry-N}, respectively, are equal up to a difference in boundary integrals in the DG Laplacian, and the location of the boundary penalty mass matrix $M_{\text{BC},h_e}$. The latter is placed in the top versus bottom diagonal block, respectively. Overall, for small $\kappa_\perp$ and $\tilde{\kappa}_{\text{BC}}$ of order $\mathcal{O}(1)$ -- as assumed in the solver analysis in \Cref{sec:solvers} -- we do not expect these two differences to affect our proposed AIR based solver's performance significantly.}{ii.}

\section{Intuition for DG-upwind operators}\label{appendix:upwind_transport}
The DG-upwind transport operator based discretization of the anisotropic heat flux can be motivated starting from an interior penalty formulation. For the anisotropic heat flux term
\begin{equation}
\nabla \cdot \big( \kappa_\Delta \mathbf{b} (\mathbf{b} \cdot \nabla T) \big),
\end{equation}
one such interior penalty formulation is given analogously to \eqref{IP_term} by
\begingroup
\addtolength{\jot}{4mm}
\begin{subequations} \label{IP_term_b}
\begin{align}
&\text{IP}_b(T_h; \phi) =- \left \langle \mathbf{b} \cdot \nabla \phi, \kappa_\Delta \mathbf{b} \cdot \nabla T_h \right \rangle \\
&+\int_\Gamma \left \llbracket (\mathbf{b} \cdot \mathbf{n}) T_h \right \rrbracket \left \lbrace \kappa_\Delta \mathbf{b} \cdot \nabla \phi \right \rbrace \;\text{d}S
+ \int_\Gamma \left \llbracket (\mathbf{b} \cdot \mathbf{n}) \phi \right \rrbracket \left \lbrace \kappa_\Delta \mathbf{b} \cdot \nabla T_h \right \rbrace \;\text{d}S \label{IP_central_b} \\
&- \int_\Gamma \frac{\kappa_\Delta \kappa_p}{h_e} \left \llbracket (\mathbf{b} \cdot \mathbf{n})\phi \right \rrbracket \left \llbracket (\mathbf{b} \cdot \mathbf{n})T_h \right \rrbracket \;\text{d}S
- \int_{\partial \Omega} \frac{2\kappa_\Delta \kappa_p}{h_e} \phi (T_h-T_{\text{BC}}) \;\text{d}S\\
&+ \int_{\partial \Omega} \kappa_\Delta (\mathbf{b} \cdot \mathbf{n}) (\mathbf{b} \cdot \nabla \phi)(T_h-T_{\text{BC}}) \;\text{d}S
+ \int_{\partial \Omega} \kappa_\Delta (\mathbf{b} \cdot \mathbf{n}) (\mathbf{b} \cdot \nabla T_h) \phi \;\text{d}S.
\end{align}
\end{subequations}
\endgroup
Next, we reformulate the central difference facet integrals \eqref{IP_central_b} to upwinded ones, further appropriately replace occurrences of $\sqrt{\kappa_\Delta} \mathbf{b}\cdot \nabla T_h$ by $\zeta_h$, and drop the interior facet penalty term. We then obtain, after some rewriting,
\begingroup
\addtolength{\jot}{4mm}
\begin{subequations} \label{IP_term_b_modified}
\begin{align}
&\widetilde{\text{IP}}_b(T_h; \phi) =- \left \langle \mathbf{b} \cdot \nabla \phi, \kappa_\Delta \mathbf{b} \cdot \nabla T_h \right \rangle - \int_{\partial \Omega} \frac{2\kappa_\Delta \kappa_p}{h_e} \phi (T_h-T_{\text{BC}}) \;\text{d}S \label{IP_term_b_kappa}\\
&+\int_\Gamma \left \llbracket \sqrt{\kappa_\Delta} T_h \mathbf{b} \cdot \mathbf{n} \right \rrbracket \sqrt{\kappa_\Delta} \widetilde{\mathbf{b} \cdot \nabla \phi } \;\text{d}S
+ \int_\Gamma \left \llbracket \sqrt{\kappa_\Delta} \phi \mathbf{b} \cdot \mathbf{n} \right \rrbracket \tilde{\zeta}_h \;\text{d}S\\
&+ \int_{\partial \Omega} \sqrt{\kappa_\Delta} (T_h-T_{\text{BC}}) (\mathbf{b} \cdot \mathbf{n}) \sqrt{\kappa_\Delta} (\mathbf{b} \cdot \nabla \phi) \;\text{d}S
+ \int_{\partial \Omega} \sqrt{\kappa_\Delta} \phi (\mathbf{b} \cdot \mathbf{n}) \zeta_h \;\text{d}S. \label{IP_term_b_modified_bdry}
\end{align}
\end{subequations}
\endgroup
Finally, we modify the boundary terms \eqref{IP_term_b_modified_bdry} in accordance with the upwind modification; to so so, we set $T_h = T_{\text{BC}}$ and $\zeta_h = \zeta_{\text{in}}$ along $\partial \Omega_{\text{in}}$. In particular, the portion along $\partial\Omega_\text{in}$ of the first integral in \eqref{IP_term_b_modified_bdry} then vanishes, and only the one along $\partial\Omega_\text{out}$ remains. Additionally, we reformulate the penalty parameter in the boundary term (last term in \ref{IP_term_b_kappa}) as $\tfrac{2\kappa_\Delta \kappa_p}{h_e} \rightarrow \kappa_\text{BC}$.

The resulting formulation can then be retrieved from our scheme \eqref{DG_upw_discr} -- ignoring the time derivative and isotropic heat flux terms -- by setting $\psi = \sqrt{\kappa}_\Delta \mathbf{b} \cdot \nabla \phi$ in \eqref{DG_upw_discr_zeta}. Note that this is merely a motivation relating the split-transport setup to an interior penalty method, since $\sqrt{\kappa}_\Delta \mathbf{b} \cdot \nabla \phi$ will generally not be in $\mathbb{V}^{\text{DG}}_k$ and we therefore cannot set $\psi$ accordingly after discretization. Importantly, however, the scheme \eqref{DG_upw_discr} still preserves the correct diffusive and total energy behavior weakly, as shown in \Cref{Remark:structure}.

\section{Base mesh}\label{appendix:base_mesh}
Here, we describe the 2D base mesh used for
\begin{figure}[ht]
\vspace{-1mm}
\begin{center}
\includegraphics[width=0.7\textwidth]{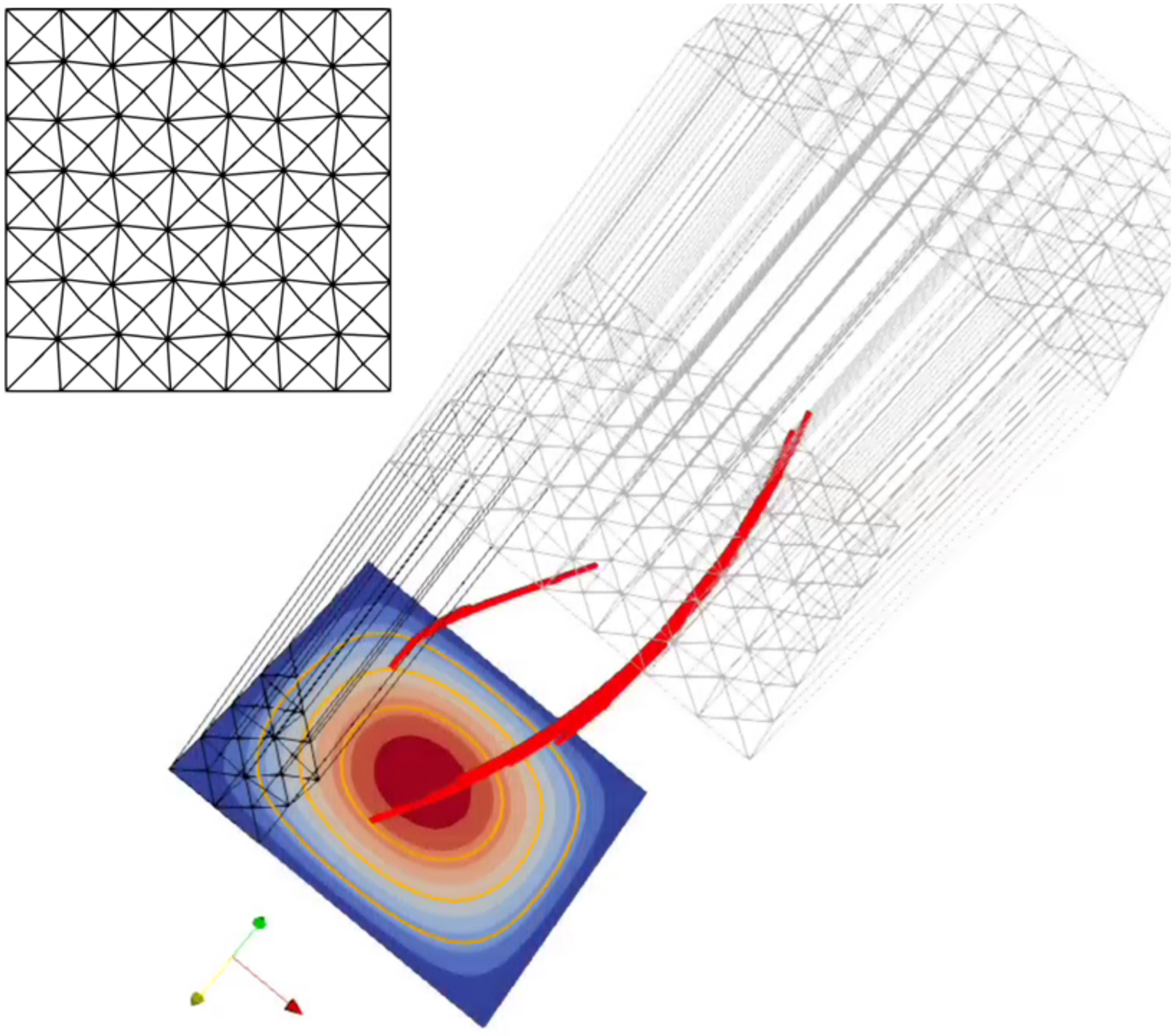}
\end{center}
\vspace{-1mm}
\caption{\ccrf{2D base mesh and extruded mesh with initial conditions \eqref{T_ic_NIMROD}, \eqref{B_ic_NIMROD} for convergence study (depicting lowest resolution). Blue to red contours on the 2D base cross-section denote the temperature field, which is constant with respect to vertical coordinate. Orange curves on the 2D base cross-section denote the magnetic field's horizontal component; red curves denote the full magnetic field's field lines.}{14}\ccbf{}{i.}} \label{fig_base_mesh_conv_ic}
\vspace{-3mm}
\end{figure}
the base square of side length $L_x$ in the numerical results section. It is generated starting from a regular triangular mesh with resolution $\Delta x = L_x/7$, whose vertices are then perturbed by a factor of $0.06\Delta x$ in each interior coordinate in order to avoid spurious error cancellations due to symmetries between the mesh and initial conditions. The resulting mesh and the corresponding initial conditions for the convergence study in \Cref{sec_conv_studies} are depicted in Figure \ref{fig_base_mesh_conv_ic}. In general, for large anisotropy ratios, we found the error values in the convergence study to be very sensitive to the underlying mesh; we therefore decided to rely on a constructed pseudo-irregular base mesh rather than a fully irregular one for better reproducibility.
%
%
\bibliographystyle{plain}
\bibliography{Aniso_open_final}
\end{document}